\documentclass[preprint,12pt]{elsarticle}

\usepackage{amssymb}
\usepackage{amsmath}

\usepackage{lineno}
\usepackage{overpic}
\journal{Applied Mathematics and Computation}

\begin{document}

\begin{frontmatter}



\title{A simple criterion for identifying slow-fast relaxation oscillations}

\author[mymainaddress]{Lucas MacQuarrie\fnref{label2}}\ead{macquarrielucas@gmail.com}

\author[mymainaddress]{Pei Yu\corref{mycorrespondingauthor}}
\ead{pyu@uwo.ca}

\affiliation[mymainaddress]{organization={Department of Mathematics, Western University},
             addressline={1151 Richmond St},
             city={London},
             postcode={N6A 5B7},
             state={Ontario},
             country={Canada}}
             
\fntext[label2]{Department of Mathematical Sciences, KAIST, Daejeon, 34101, South Korea}
\cortext[mycorrespondingauthor]{Corresponding author}




\begin{abstract}
Simple conditions have been developed in  [Zhang, Wahl and Yu, SIAM Rev. 2014; Yu and Wang, Math. Biosci. Eng. 2019], which are used to identify
the existence of slow-fast relaxation oscillations
that appear
in differential systems, where the geometric singular perturbation theory is not applicable. In this study, we present seven models which fail to satisfy these conditions, yet still show recurrence phenomena, which motivate us to generalize these conditions to establish new simple hypothetical
conditions for identifying slow-fast relaxation oscillations and providing arguments towards the development of a mathematical theory.
\end{abstract}

\begin{graphicalabstract}
\end{graphicalabstract}

\begin{highlights}
\item Research highlight 1
\item Research highlight 2
\end{highlights}

\begin{keyword}
Slow-fast motion \sep relaxation oscillation \sep recurrence \sep  
Hopf bifurcation \sep limit cycle


\MSC[2010] 34C07 \sep 34C23 \sep 34C25 \sep 34C26
\end{keyword}

\end{frontmatter}


\section{Introduction}

Systems of differential equations are able to capture a wide variety of dynamics in their solutions, making them an ideal framework for modeling real world phenomena in biology \cite{EDE2020,MUR2002,KLI2016,KAR2012}, climate \cite{LOR1963, JAC2005}, and countless other areas of science \cite{ERD1989,SUN2022,ALL1996, GOO2011, ARN1989}. In particular, systems of ordinary differential equations (ODEs), which exhibit solutions evolving over multiple timescales, are of interest as such behaviour can often be seen in models of natural phenomena. The study of such equations has inspired the development of new mathematical tools such as Geometric Singular Perturbation Theory (GSPT) \cite{Fenichel1979, JON1995, HEK2010, LI2016}. Although these techniques allow us to study a broad range of slow-fast oscillations in ODEs, there still exist systems to which the GSPT cannot be applied. 

In the following, we briefly outline the basic idea of the GSPT, which requires a given system of ODEs to be in the form of  
\begin{equation}\label{Eqn1}
\begin{array}{rl} 
\dfrac{dx}{dt} \!\!\! &= f(x,y,\varepsilon) \\[1.5ex] 
\varepsilon \dfrac{dy}{dt} &= g(x,y,\varepsilon),  
\end{array} 
\end{equation}
 where the nonlinear functions, $f\!:\mathbb{R}^n\times\mathbb{R}^m \times \mathbb{R} \to \mathbb{R}^n$ and $g\!:\mathbb{R}^n\times\mathbb{R}^m \times \mathbb{R} \to \mathbb{R}^m$, are assumed to be differentiable, and $\varepsilon$ is a small positive perturbation parameter. We see an explicit separation between the slow variable $x$ and the fast variable $y$ with appearance of $\varepsilon$. 
With a change of timescale $\tau = \varepsilon t $, the \textit{fast} system \eqref{Eqn1} becomes the \textit{slow} system,
\begin{equation}\label{Eqn1.2}
\begin{array}{rl} 
\dfrac{dx}{d\tau} &= \varepsilon f(x,y,\varepsilon) \\[1.5ex] 
 \dfrac{dy}{d\tau} &= g(x,y,\varepsilon).  
\end{array} 
\end{equation}
The systems \eqref{Eqn1} and \eqref{Eqn1.2} are equivalent for $\varepsilon >0$, but setting $\varepsilon =0$ yields the \textit{reduced equation},
\begin{equation}\label{Eqn1.3}
\begin{array}{rl} 
\dfrac{dx}{dt} \!\!\! &= f(x,y,0) \\[1.5ex] 
0 &= g(x,y,0),  
\end{array} 
\end{equation}
and the \textit{layer equation}
\begin{equation}\label{Eqn1.4}
\begin{array}{rl} 
\dfrac{dx}{d\tau} &= 0 \\[1.5ex] 
 \dfrac{dy}{d\tau} &= g(x,y,0).  
\end{array} 
\end{equation}
 Equation \eqref{Eqn1.3} is an algebraic-differential equation, describing a flow on a subset of the state space called the \textit{critical manifold} $\mathcal{M}_0$ defined by $g(x,y,0)=0$. Equation \eqref{Eqn1.4} describes the dynamics on the whole state space, however in this case the critical manifold is now comprised of equilibrium points. In other words, \eqref{Eqn1.3} describes the dynamics on the critical manifold while \eqref{Eqn1.4} describes the dynamics around it. Fenichel's theorem
\cite{Fenichel1979} tells us that if the critical manifold is both compact and normally hyperbolic, then a perturbed version of $\mathcal{M}_0$, $\mathcal{M}_\varepsilon$, exists in the system with $0<\varepsilon\ll 1$ \cite{Fenichel1979}. What's more, the stable and unstable manifolds of $\mathcal{M}_0$ persist under the perturbations as well. Thus, by studying the reduced equations, we can understand the solutions of the full system for small parameter values $\varepsilon$.

\begin{figure}[h!] 
\vspace{0.20in} 
\begin{center} 
\begin{overpic}[width=0.45\textwidth]{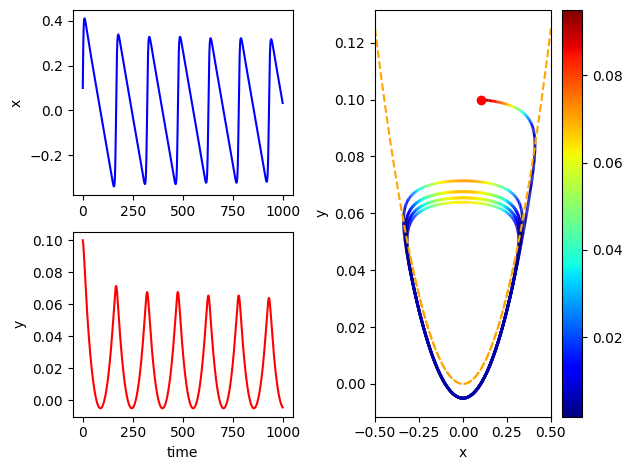}
\end{overpic} 
\hspace*{0.20in}
\begin{overpic}[width=0.45\textwidth]{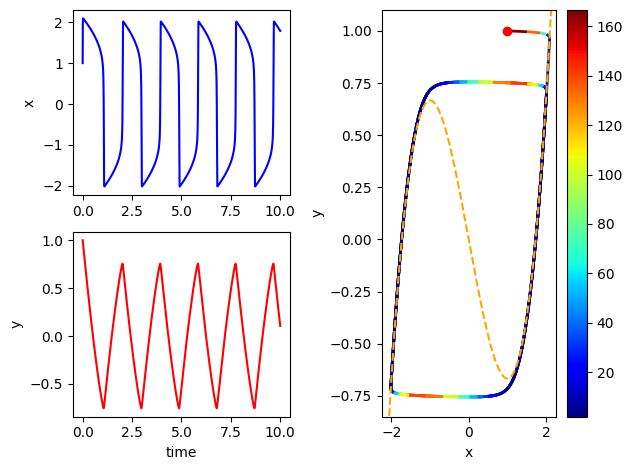}
\end{overpic} 

(a) U-shaped critical manifold  \hspace*{0.75in} (b) S-shaped critical manifold
\caption{Simulations of two systems of ODEs displaying the typical U-shape and S-shape critical manifolds as found by the GPST in dotted orange. The two ODEs are (a) $x' = y-\frac{x^2}{2}, \ y'=-\varepsilon x(1+x^2), \ \varepsilon=0.005$, constructed based on \cite{LI2013}; and (b) the Van der Pol oscillator: $x' = (x-\frac{x^3}{3}+y)/\varepsilon, \ y' = -x, \ \varepsilon=0.01$ \cite{Fenichel1979, MIS1980}. These two systems are in the form of \eqref{Eqn1}. The figures on the left side are the corresponding time histories, while the figures on the right side are the trajectories in phase space; the velocities of the solutions are indicated by solid color and the initial condition are marked with a red point.}
\label{FigPerturbed}
\end{center} 
\end{figure}

The critical manifolds of GSPT can be used to explain the dynamics of solutions of singularly perturbed equations. Figure \ref{FigPerturbed} shows two singularly perturbed systems with two typical critical manifolds. The slow-fast behaviour can be illustrated by the solutions following and leaving the critical manifold, and their difference can be seen from the shape of the critical manifold. Figure \ref{FigPerturbed} (a) shows a critical U-shaped manifold, while Figure \ref{FigPerturbed} (b) shows a critical S-shaped manifold. These shapes determine the possible behaviour of solutions as the solutions move in the state space.

It is not necessary for a system of differential equations to be of the form \eqref{Eqn1} for the
existence of slow-fast behavior.  For example, the following 2-d dimensionless HIV model \cite{ZWY2013,ZWY2014,YZ2019} 
exhibits slow-fast behaviour without being in the above form:
\begin{equation}\label{Eqn2}
\begin{array}{rl} 
\dfrac{dx}{dt} \!\!\! & = 1 - Dx - \Big(B+ \dfrac{Ay}{y+C}xy\Big), \\[1.5ex] 
\dfrac{dy}{dt} \!\!\! & = \Big(B+ \dfrac{Ay}{y+C}\Big)xy-y,
\end{array}
\end{equation}
where $x$ and $y$ represent the healthy and infected cells, respectfully, 
and $A,B,C,D$ are positive real parameters. 
\begin{figure}[h!] 
\vspace*{0.1in} 
\begin{center}
\begin{overpic}[width=0.45\textwidth]{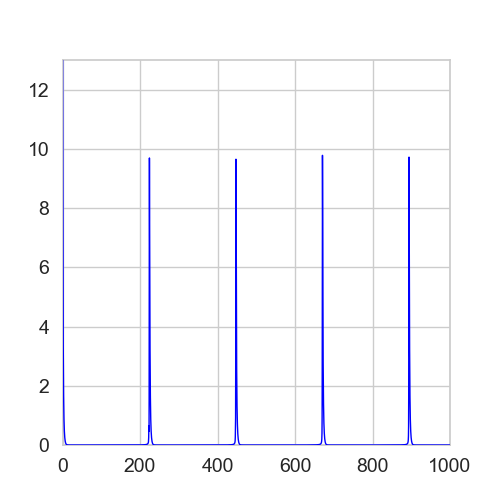}
\put(5,90){$y$} 
\put(100,6){$t$} 
\end{overpic} 
\caption{Simulated time history of model \eqref{Eqn2} for $A = 0.364$, $B=0.06$, $C = 0.823$ and $D = 0.057$ with the initial condition $(x(0),y(0))=(5,13)$, showing the slow-fast motion which is possible in a system not in the form of the singular system \eqref{Eqn1}.}
\label{addFig1}
\end{center} 
\end{figure}
The authors of \cite{ZWY2013,ZWY2014} showed that \eqref{Eqn2} exhibits a slow-fast oscillation called recurrence, as shown in Figure \ref{addFig1}, which is characterized by
long periods of quiescence followed by a quick jerk away from
equilibrium before quickly returning and then repeating. 
In the context of HIV, recurrence is the so-called viral-blips
which manifest as a spike in viral population in the patient
before returning to subcritical levels \cite{SOR2016, DIM2003}. 
Using the simple language of bifurcation theory \cite{GUC1983, WIG2003}, 
the following hypothesis with four conditions
was proposed for the existence of the recurrence phenomenon
in differential systems \cite{ZWY2013,ZWY2014}: 

\begin{itemize}
\item[{(i)}] there exist at least two equilibrium solutions;
\item[{(ii)}] there exists a transcritical bifurcation at an intersection 
of the two equilibrium solutions;
\item[{(iii)}] there is a Hopf bifurcation which occurs from one of 
the equilibrium solutions; 
\item[{(iv)}] large oscillations (or, more generally, global, 
persistent motions) can occur near the transcritical critical point.
\end{itemize}

In \cite{ZWY2013,ZWY2014}, the authors provide the explanation as follows: the slow-fast behavior occurs around a saddle point, that is, an equilibrium with one positive and one negative eigenvalue, near the transcritical bifurcation point. The saddle point connects to a stable and an unstable manifold with the corresponding eigenvalues $\lambda_1<0$ and $0<\lambda_2 \ll 1$ near the equilibrium. This indicates that the trajectories near the saddle point
will quickly approach the saddle point via the fast-stable manifold and leave slowly via the slow-unstable manifold, thus yielding slow-fast motions.  
This can also be used to explain the situation if the transcritical bifurcation condition is changed to include a saddle-node bifurcation \cite{YW2019}.

\begin{figure}[h!] 
\vspace*{0.1in} 
\begin{center}
\begin{overpic}[width=0.45\textwidth]{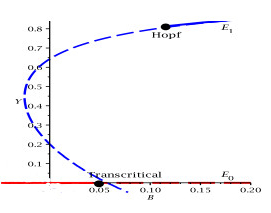}
\end{overpic} 
\caption{Bifurcation diagram of two equilibria,  $E_0$ and $E_1$, in system \eqref{Eqn2}. The branches for $E_0$ and $E_1$ are in red and blue respectively, and dashed and solid curves represent unstable adn stable regions of the equilibria respectively. Along the blue branch of $E_1$, large oscillations are allowed to exist up until the transcritical bifurcation point thus satisfying the criteria given by \cite{ZWY2013,ZWY2014}. }
\label{HIVbif}
\end{center} 
\end{figure}

Later, the four conditions were modified when the recurrent behavior was considered in a 3-d oscillating networks of biologically
relevant organic reactions \cite{YW2019,YZ2019}, which has a saddle-node bifurcation, unlike the HIV model \eqref{Eqn2} which has a transcritical bifurcation as shown in Figure \ref{HIVbif}. 

The modified 4 conditions were then proposed in \cite{YZ2019}:    
\begin{itemize}
\item[{$C_1$:}] There exists at least {\it one} equilibrium solution;
\item[{$C_2$:}] There exists a transcritical {\it or saddle-node bifurcation};
\item[{$C_3$:}] There is a Hopf bifurcation;
\item[{$C_4$:}] {\it There is a ``window'' between the Hopf bifurcation point 
and the transcritical/saddle-node bifurcation point in 
which oscillations continuously exist.} 
\end{itemize}
The main modifications include: (1) to add ``saddle-node bifurcation'' 
to the second condition; and (2) to introduce the ``window'' concept, 
which makes the last condition more clear.
The window refers to an interval $[b_T,b_H]$ or $[b_H,b_T]$ 
where $b_H$ denotes Hopf bifurcation and $b_T$ represents 
transcritical or saddle-node bifurcation, respectively. 
With the parameter values chosen inside the window, 
recurrence exists in the solutions. Thus, the application of 
the simple criterion can be applied to more systems.    

However, as we will show, there still exist systems which display the recurrence phenomenon but do not fit these modified four conditions and not all of them can be studied using the GSPT. Although the simple criterion is very powerful for its simplicity and can be easily applied to check if the recurrence phenomenon will exist with a choice of parameters, additional models in the literature have displayed recurrence behavior, 
but they do not satisfy the above conditions. 
Therefore, we are motivated to modify these conditions so
the simple criterion can be applied to more new cases, yielding 
that the previous given explanation on the sufficiency of these conditions can be adapted to our new conditions. Moreover, we present two more models that almost satisfy our conditions, resulting in recurrence-like behavior. 

In the following sections, we present model simulations
that show recurrence for a range of values of a single parameter inside a window, which are confirmed by numerically continued bifurcation diagrams \cite{KUZ2023}. 
Through these examples, we propose the following a new further modified hypotheses, which are consistent with the hypotheses $C_1$-$C_4$: 
\begin{itemize}
\item[${\rm C_1\!:}$] There exists at least one equilibrium solution;
\item[${\rm C_2\!:}$] There exists a {\bf sudden stop in oscillations};
\item[${\rm C_3\!:}$] There is a Hopf bifurcation; 
\item[${\rm C_4\!:}$] There is a ``window'' between the Hopf bifurcation point 
and the {\bf sudden stop bifurcation point} in which oscillations 
continuously exist.
\end{itemize}

\section{Five Models Displaying Recurrence without Transcritical or Saddle-node Bifurcations} 

In this section, we present five models, chosen from the current literature, which have been studied using various complex analytical methods. Here, we will confirm the existence of slow-fast oscillations by numerically simulating the models with the SciPy integrate package \cite{SCI2020}, then confirm the satisfaction of the simple criterion ${\rm C_1}$-${\rm C_4}$ through a mix of analytical and numerical methods.  For the numerical methods, our code is written in Julia \cite{JUL2017} and makes use of the BifurcationKit.jl package \cite{BIF2020}.

\subsection{ A Gause-type predator-prey system} 

\begin{figure}[h!]
\vspace*{0.1in} 
\begin{center}
\begin{overpic}[width=0.45\textwidth]{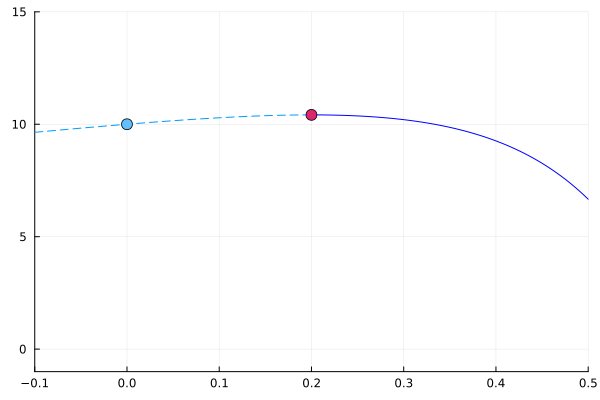}
\put(7,70){$y$} 
\put(100,4){$\varepsilon$} 
\end{overpic} 
\caption{Bifurcation diagram of the model \eqref{Eqn3} with $\varepsilon$ treated as the bifurcation parameter. The solid and dashed curves denote the stable and unstable 
equilibrium solutions, respectively.
The window can be seen between the stop point in oscillation (blue circle) and the Hopf bifurcation point (red circle).}
\label{Fig1}
\end{center} 

\end{figure}

The Gause-type predator-prey model with small death rates was considered in \cite{Hsu2019}, where the author used a modified GSPT method to prove the existence of slow-fast motion. The classical Gause-type predator-prey system with logistic growth
of the prey is described by  
\begin{equation}\label{Eqn3} 
\begin{array}{rl} 
\dfrac{dx}{dt}=\!\!\!& rx\Big(1-\dfrac{x}{K}\Big)-yp(x), \\[1.5ex] 
\dfrac{dy}{dt}=\!\!\!&  y\big[-\varepsilon +cp(x) \big],  \\[0.5ex] 
&p(0) = 0, \quad p'(0)>0, \quad p(x) > 0 \ \, \forall x>0,
\end{array}
\end{equation}
where $x$ and $y$ are the densities of prey and predator, respectively, $r$ is the intrinsic growth rate of the prey, $K$ is the carrying capacity of the environment, $c$ is the yield rate, 
and $\varepsilon$ is the death rate of the predator, all of which are positive. 
In addition, it is assumed that the function $p(x)$ is continuously
differentiable. A detailed study on the solutions of the system is given in \cite{Hsu2019}. For the purpose of simple demonstration,
take $p(x)=\frac{mx}{a+x}$ as the Holling Type-II functional response with positive constants $m$ and $a$. 
Then, the equilibrium solutions of the system are given by $(0,0)$, $(K,0)$, and 
\[(x_\epsilon, y_\epsilon)=\Big(\frac{\varepsilon a}{cm-\varepsilon},\,
\frac{r}{m}\Big(1-\frac{x_\varepsilon}{K}\Big)(x_\varepsilon+a)\Big).\] 
It is easy to see that $(x_\varepsilon, y_\varepsilon)=(0,\frac{ra}{m})$ 
 when $\varepsilon=0$. So, no transcritical bifurcation nor saddle node bifurcation can occur as long as $\frac{ra}{m}\neq K$. 
Using the BifurationKit.jl package in Julia, the branch of the bifurcation diagram for $(x_\varepsilon, y_\varepsilon)$ can be constructed numerically. We choose the parameter values
\[r=1, \ \ K=15, \ \ m=1, \ \ a=10, \ \ c=1,\]
and treat $\varepsilon$ as the bifurcation parameter. 
Examining the numerical eigenvalues evaluated at the steady-state point around $\varepsilon=0$ confirms that this point is not a transcritical bifurcation point nor a saddle-node bifurcation point, yet marks an end of the oscillations. The diagram is shown in Figure \ref{Fig1}.  

\begin{figure}[h!] 
\vspace{0.20in} 
\begin{center} 
\begin{overpic}[width=0.28\textwidth]{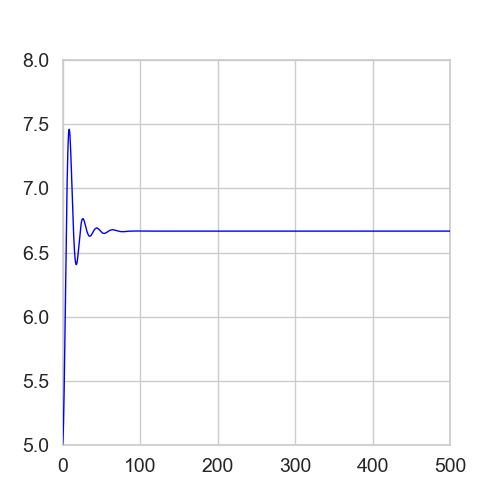} 
\put(8,95){$x$} 
\put(98,6){$t$} 
\end{overpic}
\hspace*{0.20in} 
\begin{overpic}[width=0.28\textwidth]{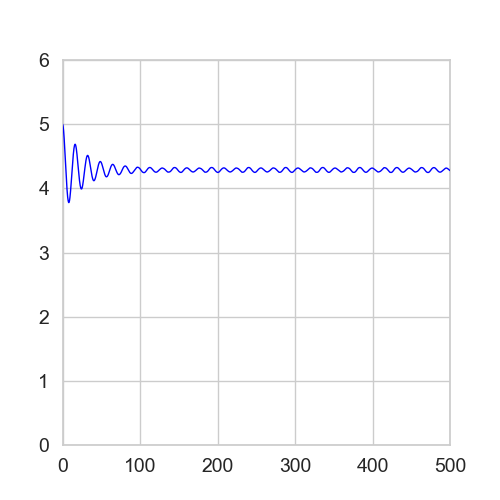}
\put(8,95){$x$} 
\put(98,6){$t$} 
\end{overpic} 
\hspace*{0.20in} 
\begin{overpic}[width=0.28\textwidth]{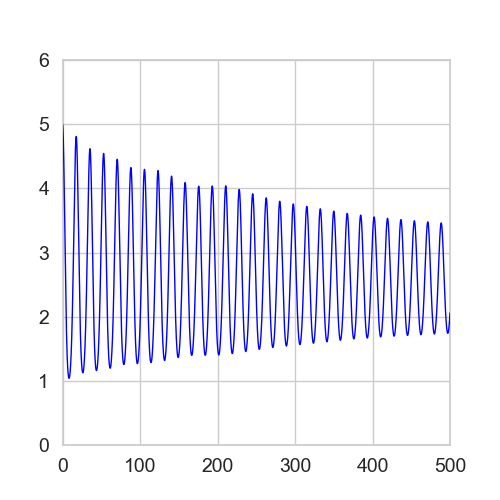}
\put(8,95){$x$} 
\put(98,6){$t$} 
\end{overpic} 

\hspace*{0.15in} (a) $\epsilon=0.4$ \hspace*{1.20in} (b) $\epsilon=0.3$
\hspace*{1.20in} (c) $\epsilon=0.2$

\vspace*{0.30in} 
\begin{overpic}[width=0.28\textwidth]{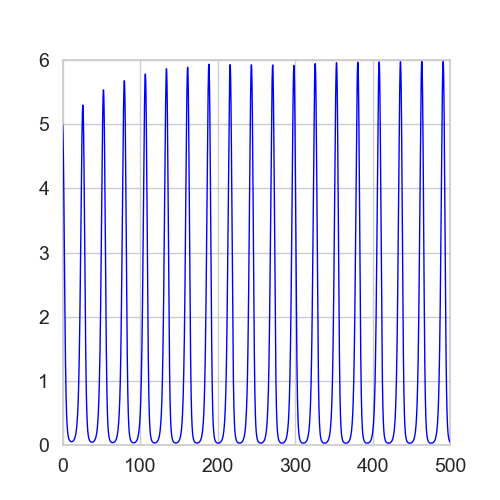} 
\put(8,95){$x$} 
\put(98,6){$t$} 
\end{overpic} 
\hspace*{0.20in} 
\begin{overpic}[width=0.28\textwidth]{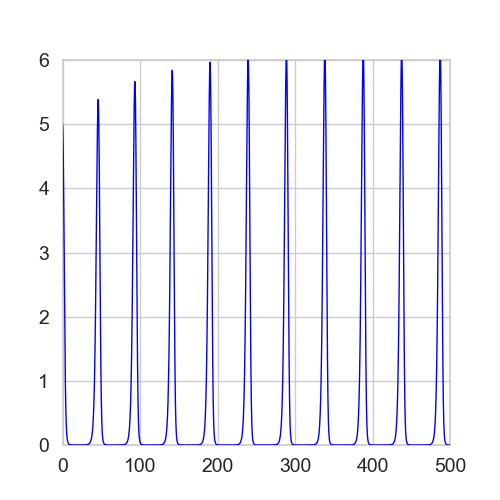} 
\put(8,95){$x$} 
\put(98,6){$t$} 
\end{overpic} 
\hspace*{0.20in} 
\begin{overpic}[width=0.28\textwidth]{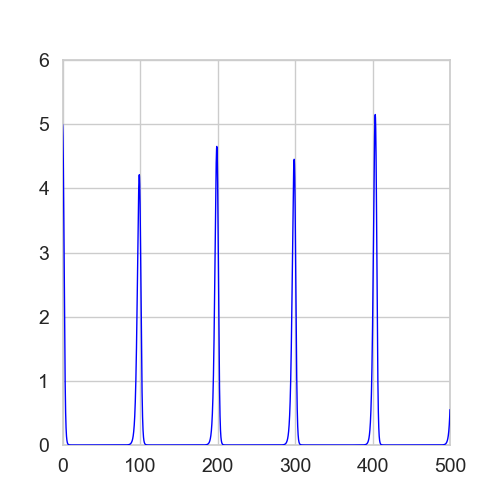} 
\put(8,95){$x$} 
\put(98,6){$t$} 
\end{overpic} 

\hspace*{0.15in} (d) $\epsilon=0.1$ \hspace*{1.20in} (e) $\epsilon=0.05$
\hspace*{1.20in} (f) $\epsilon=0.02$

\vspace*{0.56in} 
\begin{overpic}[width=0.28\textwidth]{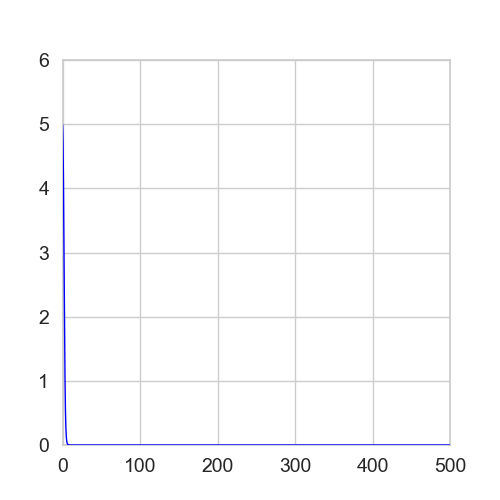} 
\put(8,95){$x$} 
\put(98,6){$t$} 
\end{overpic} 

\hspace*{0.10in} (h) $\epsilon=0$

\caption{Simulated time history of the model \ref{Eqn3} 
for $r=1$, $K=15$, $m=1$, $a=10$, $c=1$ with varying $\varepsilon$ 
and the initial point $(x(0),y(0))=(5,10)$.}
\label{Fig2}
\end{center} 
\end{figure}
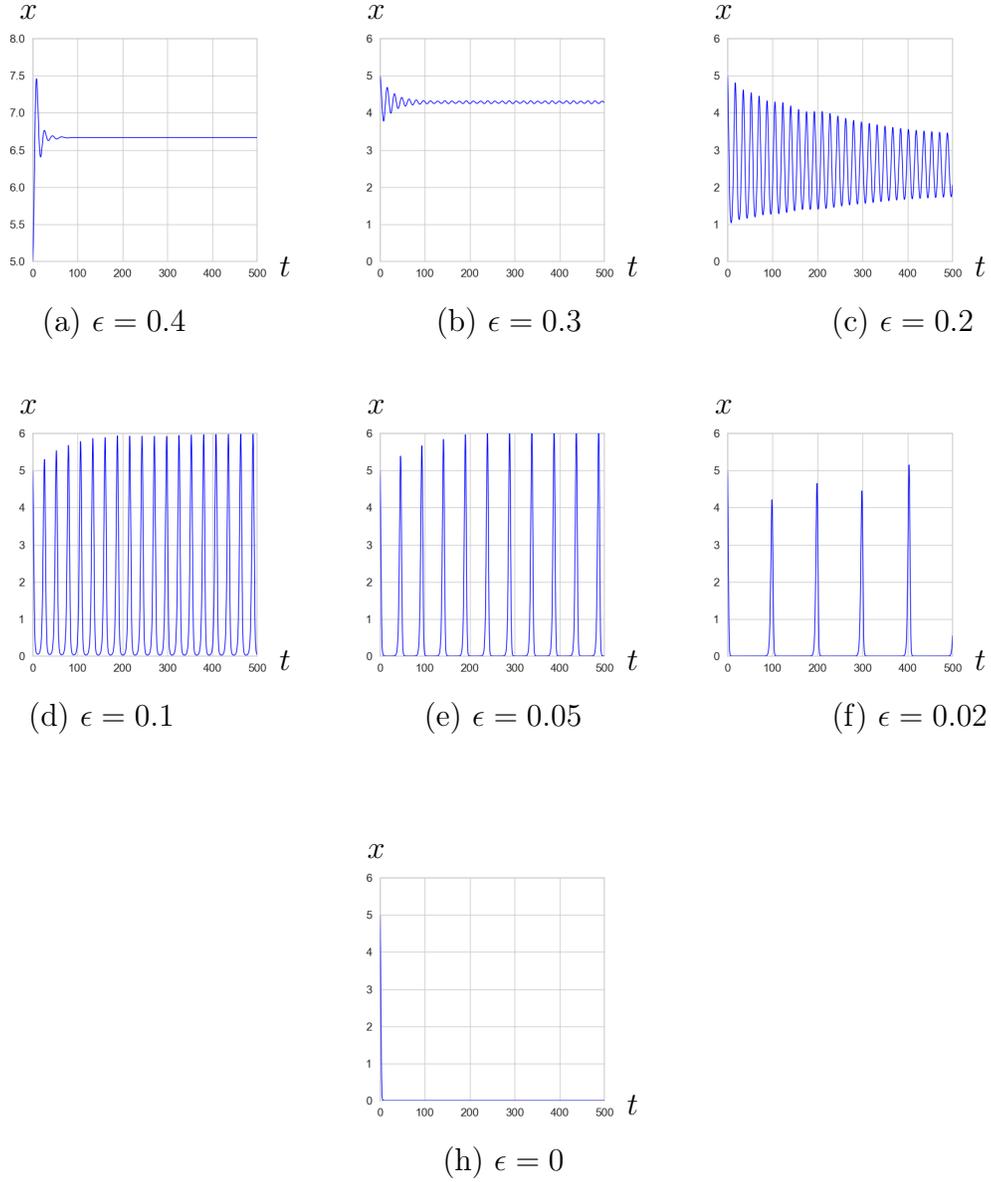

Moreover, the simulated histories for 
equation \eqref{Eqn3} are shown in figure \ref{Fig2} to confirm the existence of the slow-fast oscillation and to demonstrate the end of the motion at $\varepsilon=0$.

\subsection{A simple epidemic model}

In this subsection, we consider a simple epidemic model
\cite{Li2016}, in which the turning point plays an important role. 
The model is an SIR epidemic model incorporating demography and disease-caused death that was studied in \cite{Li2016} using the GSPT. The 3-d model is given as \cite{Li2016}
\begin{equation}\label{Eqn4} 
\begin{array}{cl} 
\dfrac{dS}{dt} & =d N+\varepsilon g(N)-h(S, N) I-(d+p) S, \\[1.0ex] 
\dfrac{dI}{dt} & = h(S, N) I-a I, \\[1.0ex] 
\dfrac{dN}{dt} & =\varepsilon g(N)-\alpha I, 
\end{array} 
\end{equation} 
where $a=d + \alpha + \gamma$ and all parameters are positive. 
It is assumed that  
\begin{enumerate}
\item $g''(N)<0$, and $g(0)=g(N^*)=0$ for some $N^*>0$. 
\item $h(S,N)$ is a smooth function and increases with 
respect to $S$ with $h(0,N)=0$. 
\item $h\big( \frac{dN}{d+p},N \big)$ is non-decreasing for 
$N\in(0,N^*)$. There exists a unique $N_0\in(0,N^*)$ such that 
$h(S_0,N_0)=a$, where $S_0=\frac{dN_0}{d+p}$ 
and $\frac{d}{d+p}h_S(S_0,N_0)+h_N(S_0,N_0)>0$. 
\end{enumerate}

For our purpose, we choose suitable functions
$h(S,N)=\frac{\beta S}{K+S}$ and $g(N)=N\left(1-\frac{N}{N_*}\right)$,
where $K$ and $N_*$ are positive parameters. 
For a bifurcation analysis, we select 
\[\alpha=0.048, \ \ \beta=1, \ \ \gamma=0.75, \ \ d=0.2, \ \ p=0.01, \ \ K=0.1, \ \ N_*=400,\]
and let $\varepsilon$ be a bifurcation parameter. 

\begin{figure}[h!] 
\vspace*{0.2in} 
\begin{center}
\begin{overpic}[width=0.45\textwidth]{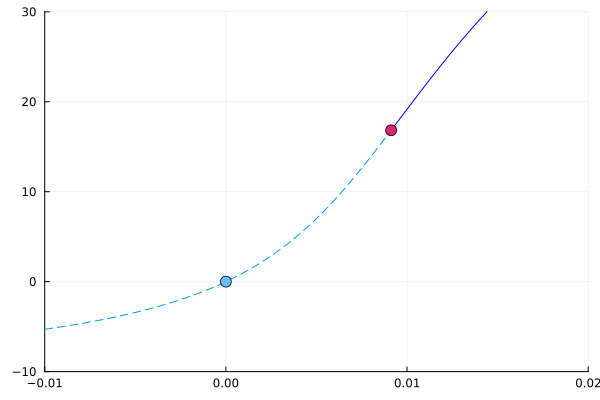}
\put(10,68){$I$} 
\put(103,4){$\varepsilon$} 
\end{overpic} 
\caption{Bifurcation diagram of the model $\eqref{Eqn4}$ projected on the $\varepsilon$-$I$
plane. The solid and dashed curves indicate the stable and unstable equilibrium solutions, respectively. The window where the slow-fast oscillation exists can be seen between the sudden stop point in the oscillation (blue circle) and the Hopf bifurcation point (red circle).}
\label{Fig3}
\end{center} 
\end{figure}

\begin{figure}[h!] 
\vspace{0.20in} 
\begin{center} 
\begin{overpic}[width=0.28\textwidth]{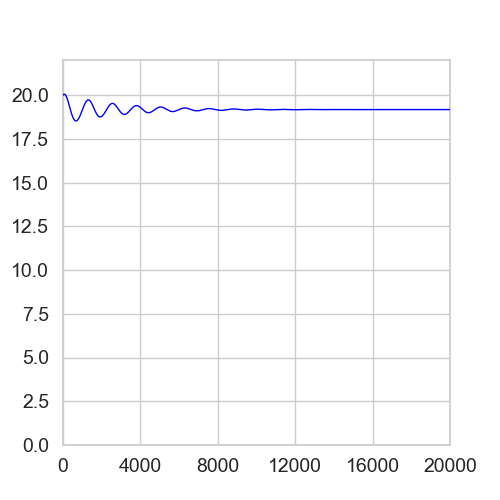} 
\put(10,90){$I$} 
\put(98,8){$t$} 
\end{overpic} 
\hspace*{0.20in} 
\begin{overpic}[width=0.28\textwidth]{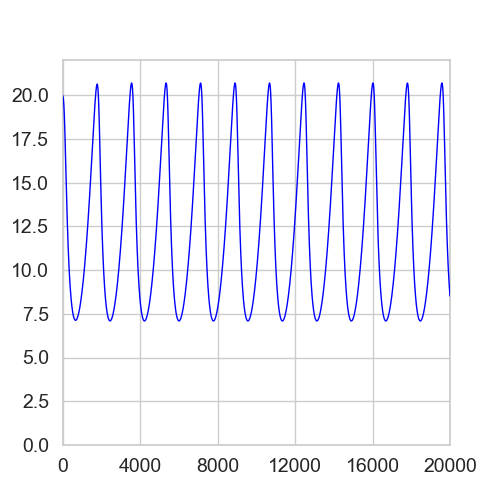}
\put(10,90){$I$} 
\put(98,8){$t$} 
\end{overpic} 
\hspace*{0.20in} 
\begin{overpic}[width=0.28\textwidth]{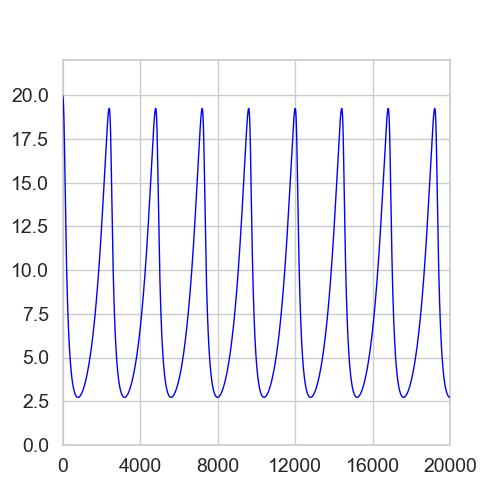}
\put(10,90){$I$} 
\put(98,8){$t$}  
\end{overpic} 

\hspace*{0.15in} (a) $\epsilon=0.01$ \hspace*{1.20in} (b) $\epsilon=0.007$
\hspace*{1.20in} (c) $\epsilon=0.005$

\vspace*{0.30in} 
\begin{overpic}[width=0.28\textwidth]{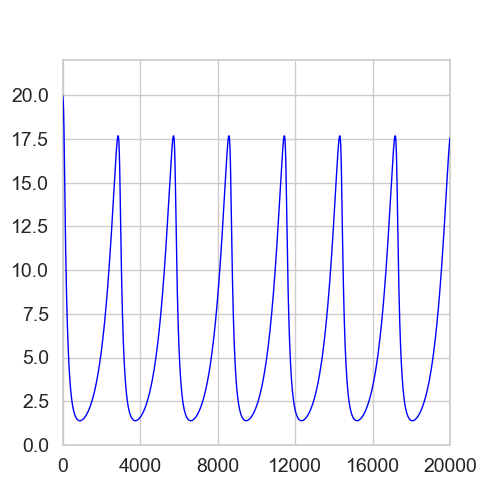}
\put(10,90){$I$} 
\put(98,8){$t$} 
\end{overpic} 
\hspace*{0.20in} 
\begin{overpic}[width=0.28\textwidth]{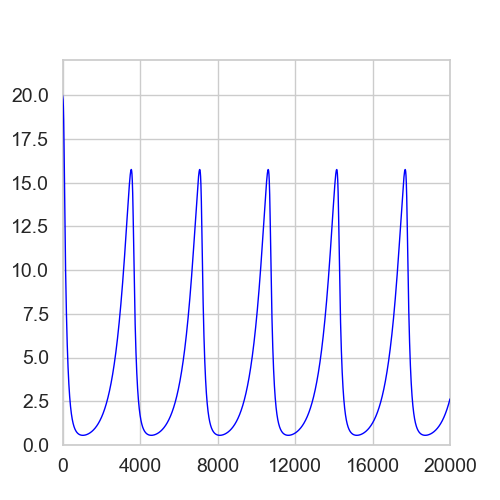}
\put(10,90){$I$} 
\put(98,8){$t$} 
\end{overpic} 
\hspace*{0.20in} 
\begin{overpic}[width=0.28\textwidth]{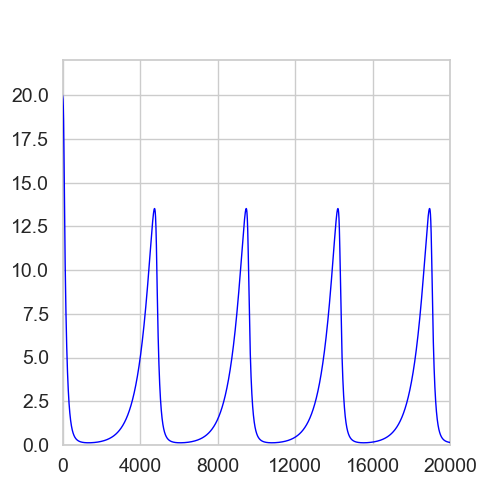}
\put(10,90){$I$} 
\put(98,8){$t$} 
\end{overpic} 

\hspace*{0.15in} (d) $\epsilon=0.004$ \hspace*{1.20in} (e) $\epsilon=0.003$
\hspace*{1.20in} (f) $\epsilon=0.002$

\vspace*{0.30in} 
\begin{overpic}[width=0.28\textwidth]{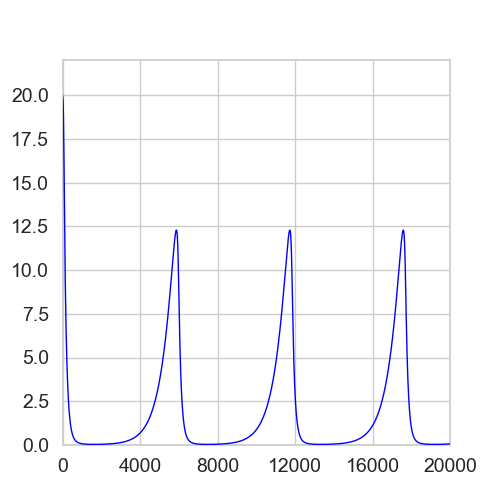}
\put(10,90){$I$} 
\put(98,8){$t$} 
\end{overpic} 
\hspace*{0.20in} 
\begin{overpic}[width=0.28\textwidth]{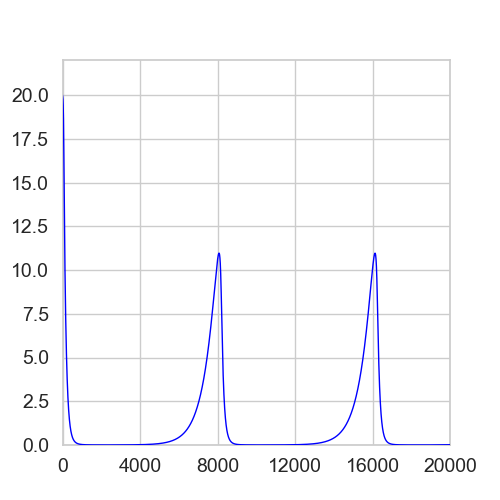}
\put(10,90){$I$} 
\put(98,8){$t$} 
\end{overpic} 
\hspace*{0.20in} 
\begin{overpic}[width=0.28\textwidth]{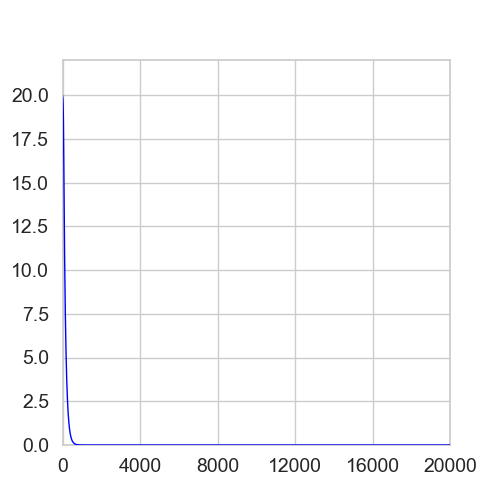}
\put(10,90){$I$} 
\put(98,8){$t$} 
\end{overpic} 
\hspace*{0.15in} (g) $\epsilon=0.0015$ \hspace*{1.20in} (h) $\epsilon=0.001$
\hspace*{1.20in} (i) $\epsilon=0$

\caption{Simulated time history of model 
\eqref{Eqn4} for $\alpha=0.048$, $\beta=1$, $\gamma=0.75$, $d=0.2$, $p=0.01$, 
$K=0.1$, and $N_*=400$, with varying $\varepsilon$ and the 
initial point: $(S(0), I(0), N(0))=(40,20,150)$.}
\label{Fig4}
\end{center} 
\end{figure}

Equation \eqref{Eqn4} is particularly interesting due to the global behavior of the system as $\varepsilon \to 0$ where the unperturbed system \eqref{Eqn4}$\mid_{\varepsilon=0}$ contains a continuum of equilibria connected by heteroclinic orbits. It has been shown in \cite{Li2016} that for $\varepsilon>0$, our choice of $h(S,N), g(N)$, and the above parameters, there exist two saddle equilibria, $(0,0,0)$ and $(S^*, 0, N^*)$, and one equilibrium $E_\varepsilon$ from which a Hopf bifurcation may occur. In particular, it is not the case that $E_\varepsilon$ approaches either of the saddle equilibria as $\varepsilon\to 0$, ensuring that this model is not covered by the conditions in \cite{YZ2019}. While \cite{Li2016} contains the conditions necessary to show the existence of the window, we numerically continue the branch of the bifurcation diagram for $E_\epsilon$ in Figure \ref{Fig3}. Moreover, we present the time history of solutions as shown in Figure \ref{Fig4} to demonstrate the recurrent behaviour. 

The slow-fast behavior of the oscillations can be clearly observed as $\varepsilon$ approaches zero. The duration of time in the state of quiescence increases as the bifurcation parameter 
$\varepsilon$ approaches the stop point in the oscillation before the oscillation is completely ceased when the parameter exits the window.

\subsection{A predator–prey model with fear and carry-over effect}

\begin{figure}
\begin{center}
\begin{overpic}[width=0.40\textwidth,height=0.20\textheight]{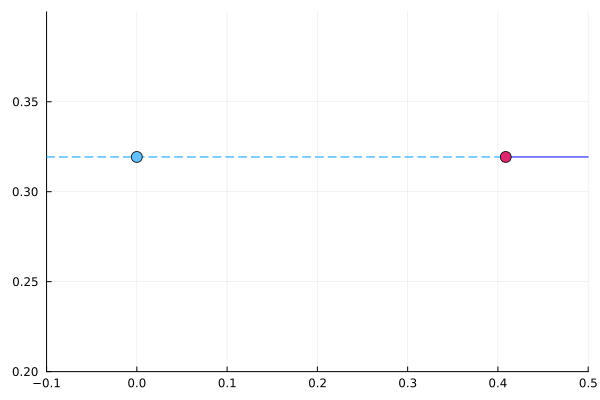}
\put(4,75){$x$} 
\put(104,5){$\varepsilon$} 
\end{overpic} 
\caption{Bifurcation diagram of the model \eqref{Eqn5} projected on the $\varepsilon$-$x$ plane. The solid and dashed lines represent stable and unstable equilibrium solutions, respectively. The window is visible between the Hopf bifurcation point (red circle) and the sudden stop point in the oscillation (blue circle).}
\label{Fig5}
\end{center}

\end{figure}

\begin{figure}[h!] 
\vspace{0.20in} 
\begin{center} 
\begin{overpic}[width=0.28\textwidth]{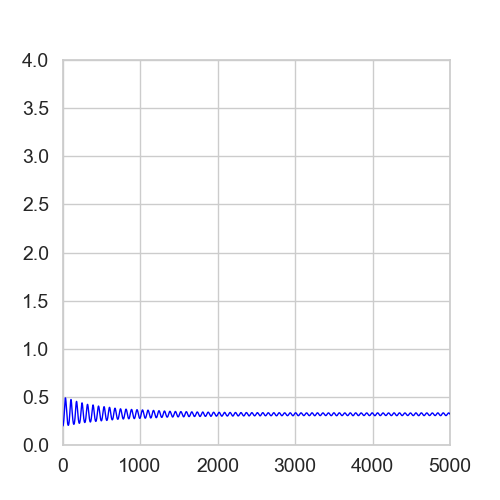} 
\put(9,93){$x$} 
\put(98,8){$t$} 
\end{overpic} 
\hspace*{0.20in} 
\begin{overpic}[width=0.28\textwidth]{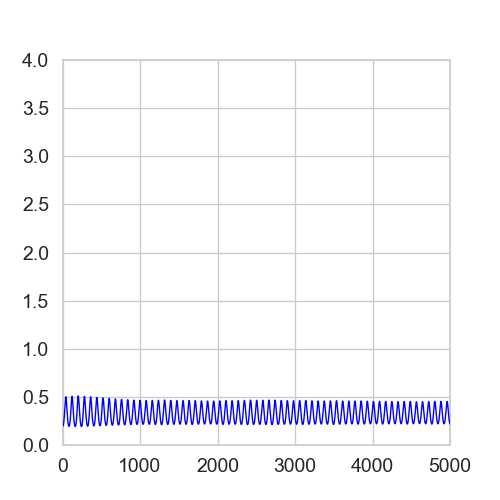} 
\put(9,93){$x$} 
\put(98,8){$t$} 
\end{overpic} 
\hspace*{0.20in} 
\begin{overpic}[width=0.28\textwidth]{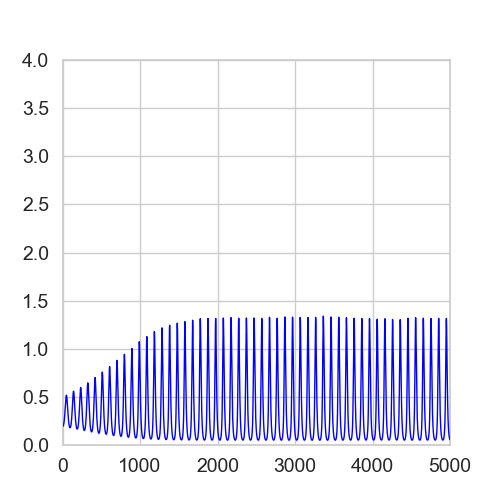} 
\put(9,93){$x$} 
\put(98,8){$t$} 
\end{overpic}
 
\hspace*{0.15in} (a) $\epsilon=0.5$ \hspace*{1.20in} (b) $\epsilon=0.4$
\hspace*{1.20in} (c) $\epsilon=0.3$

\vspace*{0.30in} 
\begin{overpic}[width=0.28\textwidth]{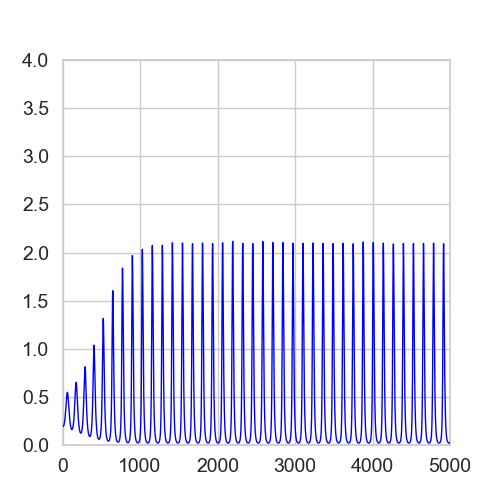} 
\put(9,93){$x$} 
\put(98,8){$t$} 
\end{overpic} 
\hspace*{0.20in} 
\begin{overpic}[width=0.28\textwidth]{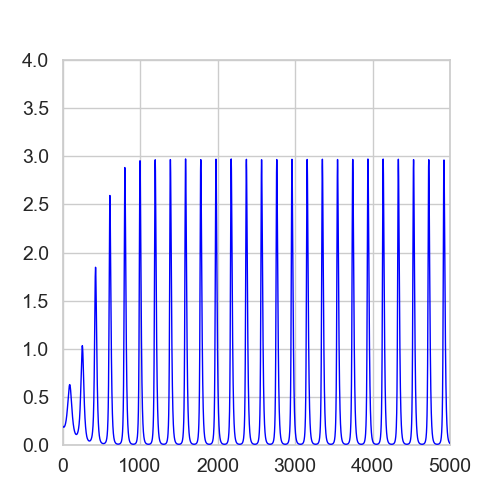} 
\put(9,93){$x$} 
\put(98,8){$t$} 
\end{overpic} 
\hspace*{0.20in} 
\begin{overpic}[width=0.28\textwidth]{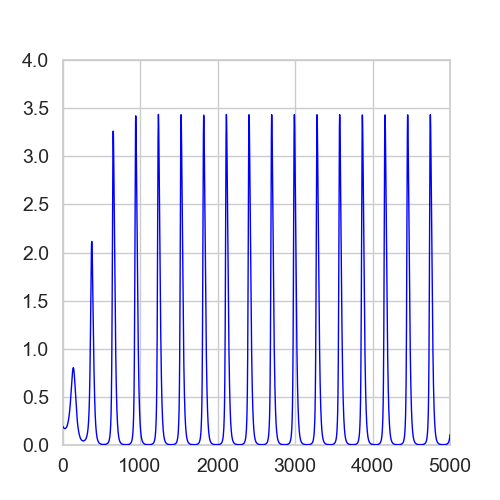} 
\put(9,93){$x$} 
\put(98,8){$t$} 
\end{overpic} 

\hspace*{0.15in} (d) $\epsilon=0.2$ \hspace*{1.20in} (e) $\epsilon=0.1$
\hspace*{1.20in} (f) $\epsilon=0.05$

\vspace*{0.30in} 
\begin{overpic}[width=0.28\textwidth]{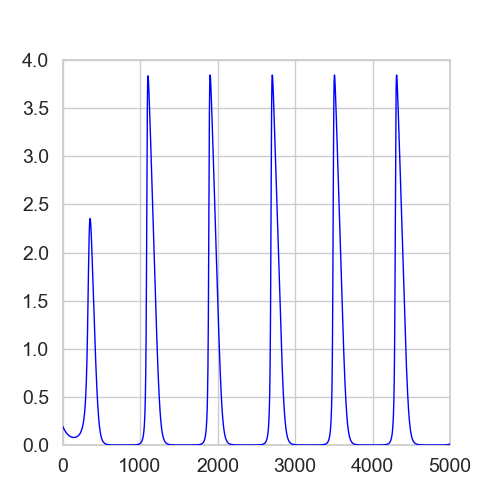} 
\put(9,93){$x$} 
\put(98,8){$t$} 
\end{overpic} 
\hspace*{0.20in} 
\begin{overpic}[width=0.28\textwidth]{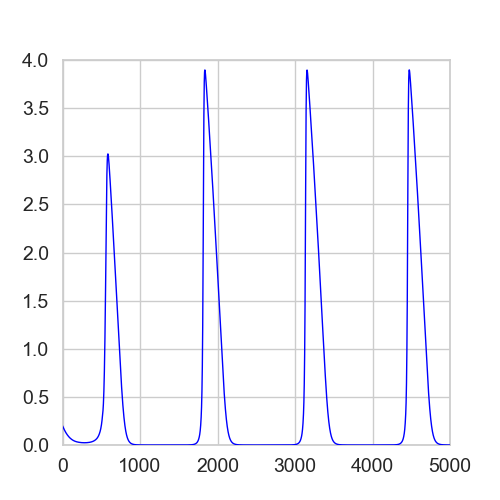} 
\put(9,93){$x$} 
\put(98,8){$t$}  
\end{overpic} 
\hspace*{0.20in} 
\begin{overpic}[width=0.28\textwidth]{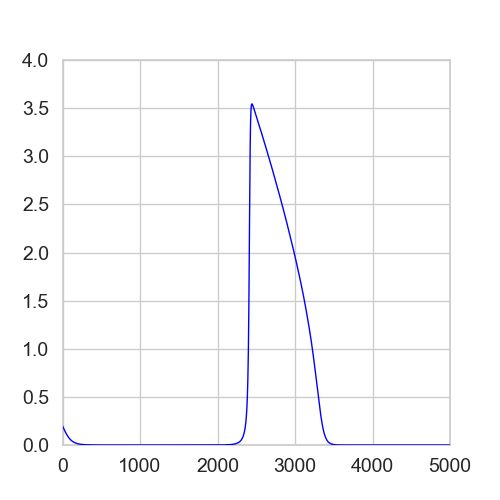} 
\put(9,93){$x$} 
\put(98,8){$t$} 
\end{overpic} 

\hspace*{0.15in} (g) $\epsilon=0.01$ \hspace*{1.20in} (h) $\epsilon=0.005$
\hspace*{1.20in} (i) $\epsilon=0.001$

\caption{Simulated time histories of the model \eqref{Eqn5} under the variation of $\epsilon$, with the initial point $(x(0),y(0))=(0.2,0.9).$}
\label{Fig6}
\end{center} 
\end{figure}

A predator-prey model with fear and carry-over effect showing
slow-fast oscillations is described by the following 
dimensionless differential equations in \cite{Sahoo-Samanta2023}: 
\begin{equation}\label{Eqn5} 
\begin{array}{rl} 
\dfrac{dx}{dt} \!\!\! &=x \Big[ \dfrac{1+x}{1+x+\kappa y}-\delta_1 
-\gamma x \Big], \\[1.5ex] 
\dfrac{dy}{dt} \!\!\! &=\varepsilon y\Big(\dfrac{x}{\theta+y+x}-\delta_2\Big).
\end{array}
\end{equation} 
The parameter values selected in \cite{Sahoo-Samanta2023} for simulations are
\[\delta_1=0.2, \ \ \delta_2=0.1, \ \ \kappa=5, \ \ \gamma=0.1, \ \  \theta=2, \]
and $\varepsilon$ is treated as a bifurcation parameter. 
It was shown in \cite{Sahoo-Samanta2023} that three equilibria can exist: 
a trivial equilibrium $E_0 = (0,0)$, a boundary equilibrium 
$E_a = (\frac{1-\delta_1}{\gamma},0)$ and 
a unique coexistence equilibrium $E_*=(x^*,y^*)$, where $y_*$ is a non-zero component. All equilibria are independent of $\epsilon$ (if they exist). 
Thus, neither a transcritical bifurcation nor a saddle-node bifurcation can exist. 
As discussed in the previous sections, the bifurcation diagram was obtained using a numerical integration
scheme to show the existence of the theorized window between the branch
point and the Hopf bifurcation point, as shown in Figure \ref{Fig5}. 
Again, simulated time histories of this model with varying $\varepsilon$ are given in Figure \ref{Fig6}, exhibiting the recurrence behavior.

\subsection{A three-species Lotka–Volterra food web model}

\begin{figure}[!h] 
\vspace*{0.5in} 
\begin{center}
\begin{overpic}[width=0.45\textwidth]{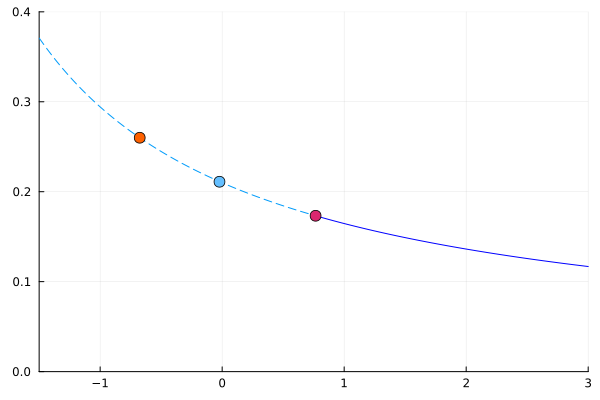} 
\put(8,70){$y$} 
\put(103,4){$\beta$} 
\end{overpic} 
\hspace*{0.20in} 
\begin{overpic}[width=0.45\textwidth]{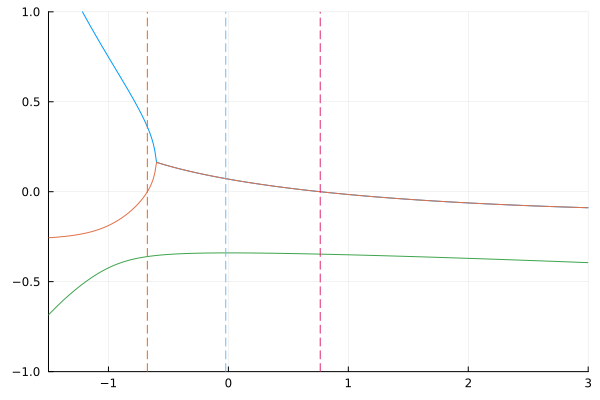}
\put(7,70){$\lambda_i$} 
\put(103,4){$\beta$} 
\end{overpic} 

(a) Bifurcation diagram of the model \eqref{Eqn6} \hspace*{0.40in} (b) Real part of corresponding eigenvalues to (a).
\caption{(a) Bifurcation diagram of the model \eqref{Eqn6} for one of the equilibria, projected on the $\beta$-$y$ plane; and (b) the real part of the three eigenvalues $\lambda_i$, ($i=1,2,3$), (solid lines) corresponding to the
equilibrium solution in part (a). The solid and dashed curves in part (a) denote the stable and unstable equilibrium solutions, respectively. The window for the slow-fast oscillation is visible between the sudden stop point in the oscillation (blue circle) and the Hopf bifurcation point (red circle). A branch point (orange circle) is shown for completeness of the branch and does not effect the window. In part (b), the dashed verticle lines indicate the bifurcation points of the equilibrium, where the colors match their respective points in part (a). The color of the solid lines in part (b) distinguish each of the eigenvalues of the Jacobian of the system evaluated at the equilibrium point. It is seen that the sudden stop point in the oscillation is not associated with a change in the sign of the real part of any eigenvalues, implying that the sudden stop of oscillation is not due to a transcritical or a saddle-node bifurcation.}
\label{Fig7}
\end{center}

\end{figure}

We now consider a three-species Lotka–Volterra food web model 
with omnivory studied in \cite{Hsu2015}. The non-dimensionalized equations are given by
\begin{equation}\label{Eqn6}
\begin{array}{rl} 
\dfrac{dx}{dt} \!\!\! &= x(1-x-y-\bar{\gamma}z), \\[1.0ex] 
\dfrac{dy}{dt} \!\!\! &= y(-d_1+\alpha x - \beta z), \\[1.0ex]
\dfrac{dz}{dt} \!\!\! &= z(-d_2 + \gamma x + \delta y),  
\end{array} 
\end{equation}
where the initial conditions and the parameters 
$\bar{\gamma}$, $\alpha$, $\beta$, $\gamma$, $\delta$, $d_1$ and $d_2$ 
are all positive. Fixing 
$$
\alpha=2.5, \ \ \bar{\gamma}=1, \ \ \gamma=0.25, \ \ 
\delta = 1, \ \ d_1=0.5, \ \ d_2=0.26, 
$$
and choosing $\beta$ as the bifurcation parameter, we obtain a simulated bifurcation diagram for the positive equilibrium from which the limit cycle bifurcates, as shown in Figure \ref{Fig7}(a).   
While $\beta$ is specified to be positive for ecological reasons, we let $\beta$ become negative which allows us to conclude that
the recurrence phenomenon  correlates with a sudden stop in the oscillation
around $\beta \approx -0.022$. As discussed in the previous models, this leads to a new type of window in the bifurcation diagram after the Hopf bifurcation, which was not covered in \cite{YZ2019}. Numerically calculating the real part of the eigenvalues (depicted in \ref{Fig7}(b)) shows that none of them are zero around $\beta \approx -0.022$, confirms that the stop in the oscillation is not due to a saddle-node or transcritical bifurcation. 

\begin{figure}[h!] 
\vspace{0.20in} 
\begin{center} 
\begin{overpic}[width=0.28\textwidth]{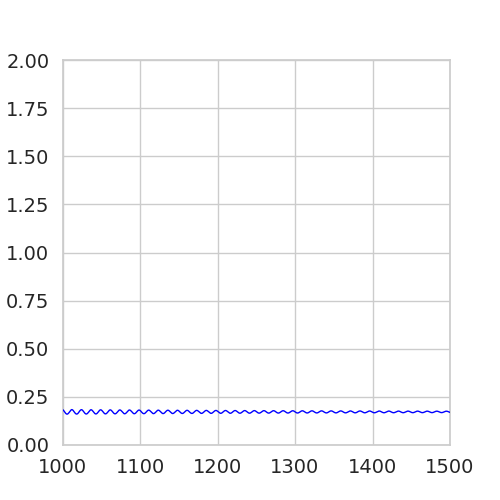} 
\put(10,95){$y$} 
\put(98,6){$t$} 
\end{overpic}
\hspace*{0.20in} 
\begin{overpic}[width=0.28\textwidth]{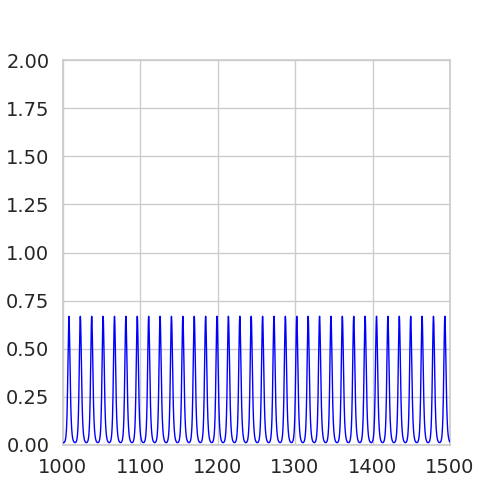}
\put(10,95){$y$} 
\put(98,6){$t$} 
\end{overpic} 
\hspace*{0.20in} 
\begin{overpic}[width=0.28\textwidth]{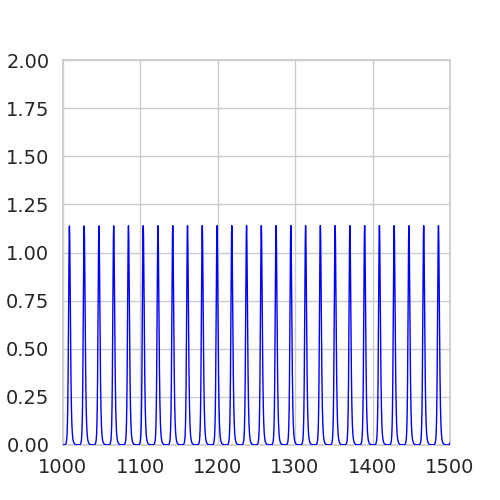} 
\put(10,95){$y$} 
\put(98,6){$t$} 
\end{overpic} 

\hspace*{0.15in} (a) $\beta=0.8$ \hspace*{1.20in} (b) $\beta=0.6$
\hspace*{1.20in} (c) $\beta=0.4$

\vspace*{0.30in} 
\begin{overpic}[width=0.28\textwidth]{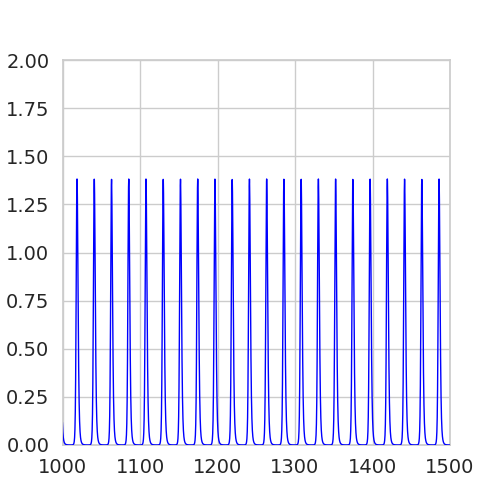} 
\put(10,95){$y$} 
\put(98,6){$t$} 
\end{overpic} 
\hspace*{0.20in} 
\begin{overpic}[width=0.28\textwidth]{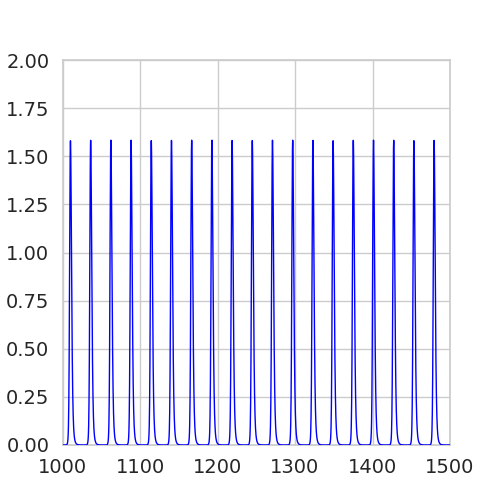}
\put(10,95){$y$} 
\put(98,6){$t$} 
\end{overpic} 
\hspace*{0.20in} 
\begin{overpic}[width=0.28\textwidth]{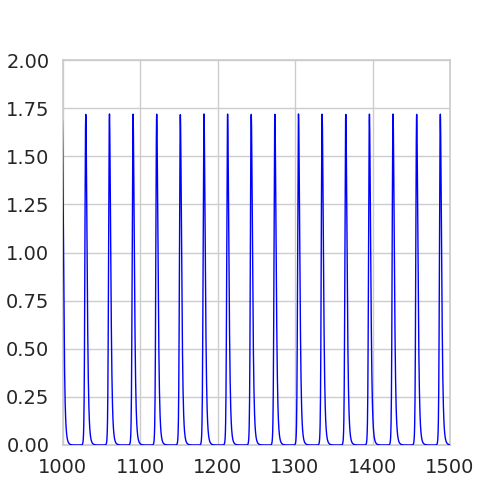} 
\put(10,95){$y$} 
\put(98,6){$t$} 
\end{overpic} 

\hspace*{0.15in} (d) $\beta=0.3$ \hspace*{1.20in} (e) $\beta=0.2$
\hspace*{1.20in} (f) $\beta=0.1$

\vspace*{0.30in} 
\begin{overpic}[width=0.28\textwidth]{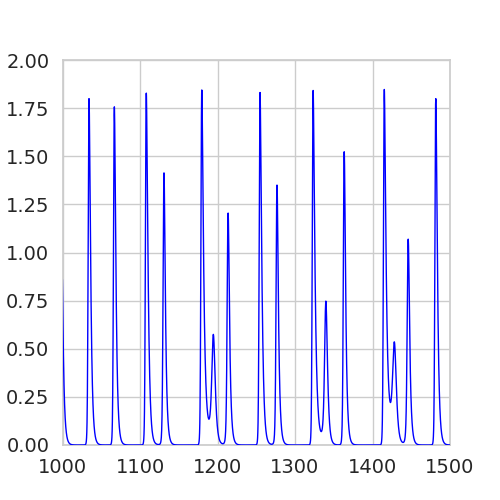}
\put(10,95){$y$} 
\put(98,6){$t$} 
\end{overpic} 
\hspace*{0.20in} 
\begin{overpic}[width=0.28\textwidth]{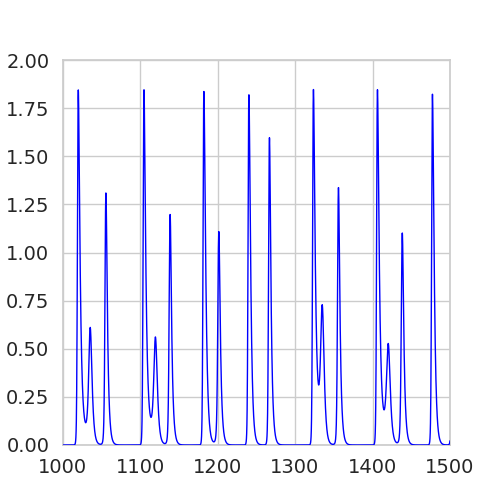} 
\put(10,95){$y$} 
\put(98,6){$t$} 
\end{overpic} 
\hspace*{0.20in} 
\begin{overpic}[width=0.28\textwidth]{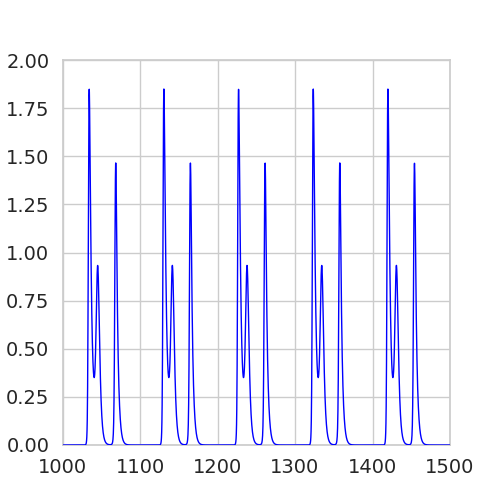}
\put(10,95){$y$} 
\put(98,6){$t$} 
\end{overpic} 

\hspace*{0.15in} (g) $\beta=0.01$ \hspace*{1.20in} (h) $\beta=0.001$
\hspace*{1.20in} (i) $\beta=0$

\caption{Simulated time histories of the model \eqref{Eqn6} for $\alpha=2.5$, $\bar{\gamma}=1$, $\gamma=0.25$, $\delta = 1$, 
$d_1=0.5$ and $d_2=0.26$, with varying $\beta$ and the 
initial point: $(x(0),y(0),z(0))=(1,1,1)$, showing the birth of the limit cycles (see Figure \ref{Fig7}), followed by progressively more pronounced slow-fast oscillations near the sudden stop point in the oscillation.}
\label{Fig8}
\end{center} 

\end{figure}
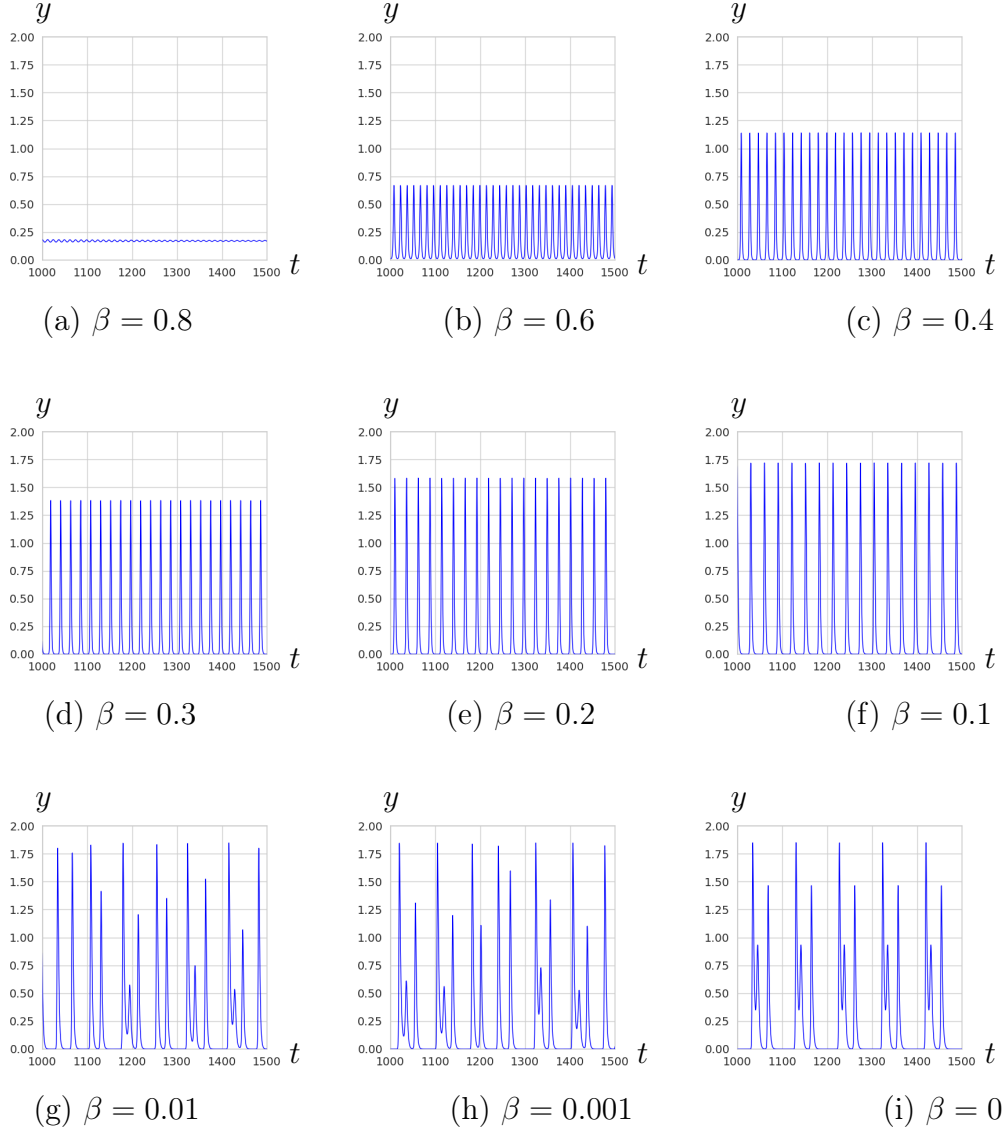

Simulating the system for $\beta$ slightly greater than $-0.022$ shows the existence of the slow-fast oscillation which is unique among those models we have studied thus far. The various peaks in Figures \ref{Fig8}(g), (h), and (i), 
hint at the existence of a limit cycle with many loops, indicating that our hypotheses apply to more than the simple periodic motion we have seen in Sections 2.1, 2.2, and 2.3. One can also notice that the period of the oscillations does not seem to increase as the bifurcation parameter approaches the stop point, which was the case shown in Figures \ref{Fig4} and \ref{Fig6}. 

\subsection{ An El Ni$\tilde{\rm n}$o Southern oscillation model}

In this subsection, we consider complex oscillatory dynamics in a $3$-timescale El Ni$\tilde{\rm n}$o Southern Oscillation (ENSO) model
\cite{Kaklamanos2023}, which describes the global climate dynamics. 
The ENSO model is given by the following dimensionless differential 
equations:
\begin{equation}\label{Eqn7}
\begin{array}{ll} 
\dfrac{dx}{dt}=x[x+y+c(1-\tanh \{x+z\})]+\rho \delta\left(x^2-a x\right),
\\[1.0ex] 
\dfrac{dy}{dt}=-\rho \delta\left(a y+x^2\right), \\[1.0ex]
\dfrac{dz}{dt}=\delta\left(k-z-\frac{x}{2}\right),
\end{array}
\end{equation} 
where the state variables and parameters are defined as
$$
\begin{aligned}
& x \leq 0, \quad y \in \mathbb{R}, \quad z \geq 0, \\
& c \in\left(1, c_0\right), \quad k \in(0,1), \quad a \in\left(0, a_0\right), 
\quad \text { and } \quad 0<\delta, \rho \ll 1,
\end{aligned}
$$
for some fixed $c_0>1$ and $a_0>0$.
As explained in \cite{Kaklamanos2023},  
$x$ corresponds to the temperature difference between the Eastern 
and Western Pacific surface water; $y$ corresponds to the departure temperature of the Western Pacific surface ocean from certain reference mean temperature; and $z$ represents the Western Pacific thermocline depth anomaly. Therefore, it is apparent that the parameters: $c$, $k$, $a$, $\delta$ and $\rho$ are associated with the rates of the aforementioned processes.

In the following, we will show that the sudden stop in oscillations cannot be attributed to a transcritical or saddle-node bifurcation of an equilibrium and that the newly hypothesized window exists. To illustrate this, we fix 
\[a=2, \ \ c=1.4, \ \ k=0.7, \ \ \rho=0.01, \]
and take $\delta$ as a bifurcation parameter
to demonstrate that oscillations cannot exist when $\delta=0$. In this case, $\frac{dy}{dt}=\frac{dz}{dt}=0$ and the non-trivial dynamics are reduced to a first-order nonlinear ODE: $$
\dfrac{dx}{dt} =x(x+y(0)+c(1-\tanh(x-z(0))),
$$
whose solution cannot be periodic, implying that the oscillation stops at $\delta=0$. 

We can find the equilibria of this system by setting the derivatives to zero. The equations for $\frac{dy}{dt}$ and $\frac{dz}{dt}$ imply $y = -\frac{x^2}{a}$ and $z = k-\frac{x}{2}$, respectively, which are then substituted into the equation for $\frac{dx}{dt}$
to yield the following equation after simplification, 
\begin{equation}
 x\left[ -\frac{x^2}{a} + (1+\rho \delta)x -\rho \delta a +c\left(1-\tanh\left(\frac{x}{2}+k\right)\right)\right] = 0. \label{equilibrium}
\end{equation}
For our chosen parameter values and with $\delta$ sufficiently small, \eqref{equilibrium} has two solutions satisfying $x \leq 0 $: $x=0$ and $x=x_0$. The equilibria are thus given by $E_0=(0,0,k)$ and $E_1=\left(x_0, -\frac{x_0^2}{a}, k-\frac{x_0}{2})\right)$, where the value of $x_0$ depends upon the parameters. To ensure that a transcritical or saddle-node bifurcation does not exist, one must have that these two equilibria do not coincide at $\delta=0$. Setting $\delta =0$ in \eqref{equilibrium} yields $-\frac{x_0^2}{a} + x +c(1-\tanh (\frac{x_0}{2}+k))=0$ whose solution is not $x_0=0$ and therefore $E_0 \neq E_1$. Thus, $\delta_0$ is not a transcritical bifurcation point nor a saddle-node bifurcation point.

\begin{figure}[h!] 
\vspace*{0.2in} 
\begin{center}
\begin{overpic}[width=0.45\textwidth]{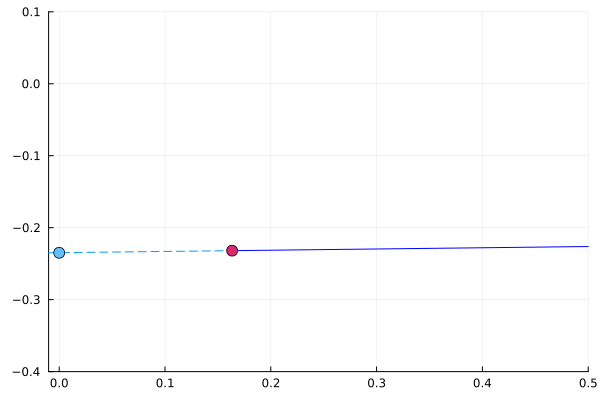}
\put(19,73){$y$} 
\put(100,5){$\delta$} 
\end{overpic} 
\caption{Bifurcation diagram of the model \eqref{Eqn7} projected on the $\delta$-$y$ plane. The solid and dashed lines denote the stable and unstable equilibrium solutions, respectively. The window for the slow-fast oscillation is between the sudden stop point in the  oscillation (blue circle) and the Hopf bifurcation point (red circle).}
\label{Fig9}
\end{center}

\end{figure}

\begin{figure}[h!] 
\vspace{0.20in} 
\begin{center} 
\begin{overpic}[width=0.28\textwidth]{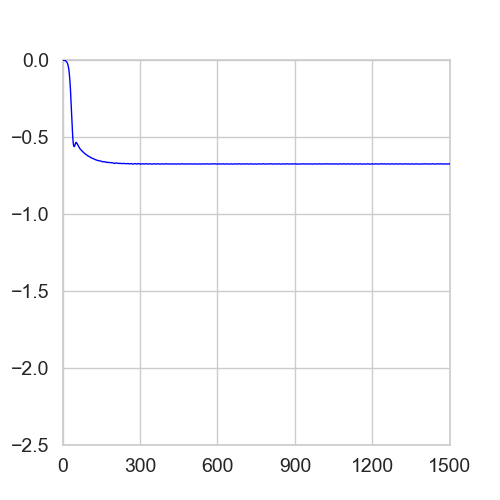} 
\put(10,90){$x$} 
\put(98,8){$t$} 
\end{overpic}
\hspace*{0.20in} 
\begin{overpic}[width=0.28\textwidth]{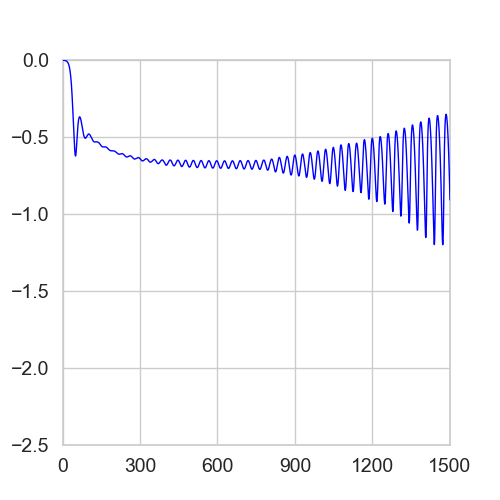}
\put(10,90){$x$} 
\put(98,8){$t$} 
\end{overpic} 
\hspace*{0.20in} 
\begin{overpic}[width=0.28\textwidth]{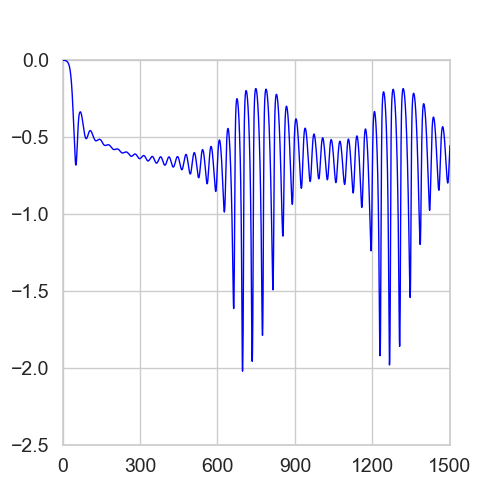}
\put(10,90){$x$} 
\put(98,8){$t$} 
\end{overpic} 

\hspace*{0.15in} (a) $\delta=0.4$ \hspace*{1.20in} (b) $\delta=0.16$
\hspace*{1.20in} (c) $\delta=0.14$

\vspace*{0.30in} 
\begin{overpic}[width=0.28\textwidth]{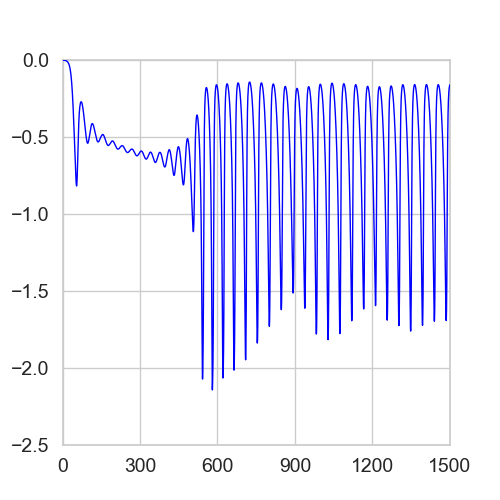} 
\put(10,90){$x$} 
\put(98,8){$t$} 
\end{overpic} 
\hspace*{0.20in} 
\begin{overpic}[width=0.28\textwidth]{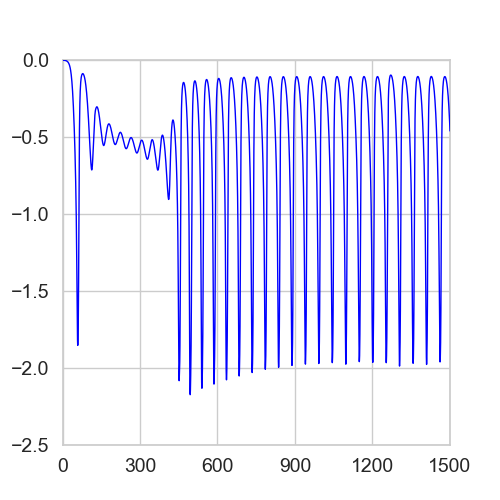} 
\put(10,90){$x$} 
\put(98,8){$t$} 
\end{overpic} 
\hspace*{0.20in} 
\begin{overpic}[width=0.28\textwidth]{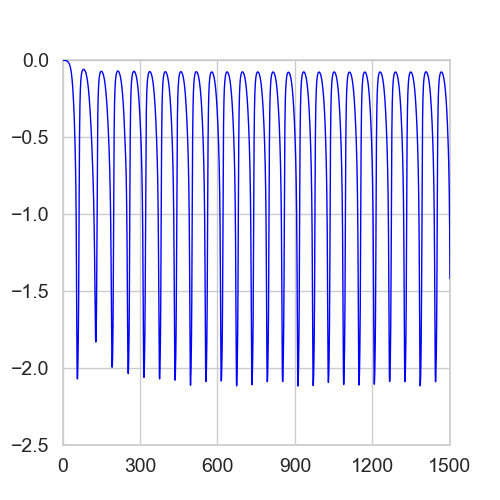} 
\put(10,90){$x$} 
\put(98,8){$t$} 
\end{overpic} 

\hspace*{0.15in} (d) $\delta=0.12$ \hspace*{1.20in} (e) $\delta=0.1$
\hspace*{1.20in} (f) $\delta=0.08$

\vspace*{0.30in} 
\begin{overpic}[width=0.28\textwidth]{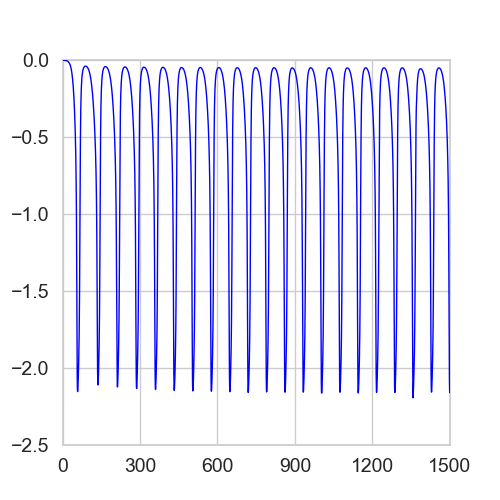} 
\put(10,90){$x$} 
\put(98,8){$t$} 
\end{overpic} 
\hspace*{0.20in} 
\begin{overpic}[width=0.28\textwidth]{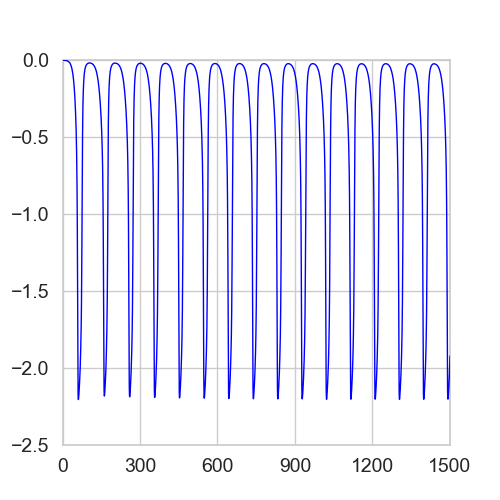} 
\put(10,90){$x$} 
\put(98,8){$t$}  
\end{overpic} 
\hspace*{0.20in} 
\begin{overpic}[width=0.28\textwidth]{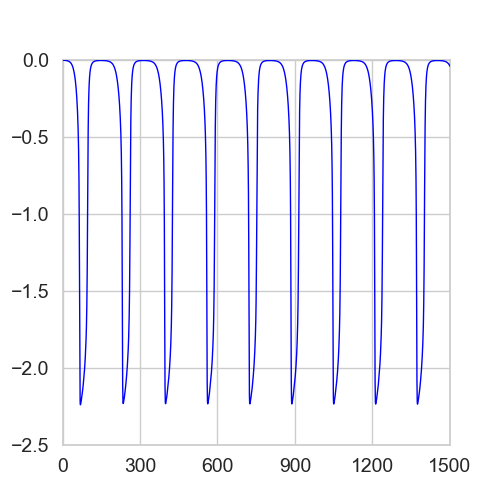} 
\put(10,90){$x$} 
\put(98,8){$t$} 
\end{overpic} 

\hspace*{0.15in} (g) $\delta=0.06$ \hspace*{1.20in} (h) $\delta=0.04$
\hspace*{1.20in} (i) $\delta=0.02$

\vspace*{0.30in} 
\begin{overpic}[width=0.28\textwidth]{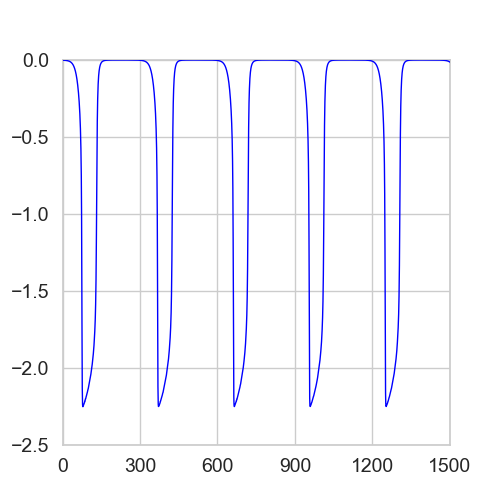} 
\put(10,90){$x$} 
\put(98,8){$t$} 
\end{overpic} 
\hspace*{0.20in} 
\begin{overpic}[width=0.28\textwidth]{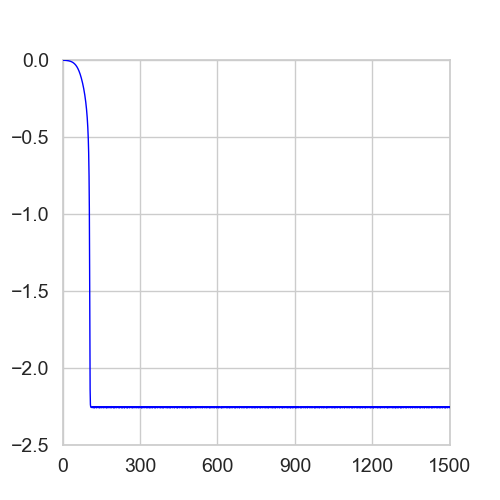} 
\put(10,90){$x$} 
\put(98,8){$t$} 
\end{overpic} 

\hspace*{0.15in} (j) $\delta=0.01$ \hspace*{1.20in} (k) $\delta=0$

\caption{Simulated time histories of the model \eqref{Eqn7} under the variation of $\delta$, with the initial point $(x(0),y(0),z(0))=(-0.001,-0.4,0.8)$}
\label{Fig10}
\end{center} 

\end{figure}

To this end, we shall use a numerical bifurcation analysis as shown in Figure \ref{Fig9} to demonstrate the existence of a Hopf bifurcation from $E_1$ at $\delta \approx 0.1634$. The simulated time histories are shown in Figure \ref{Fig10} from which one can see the periods of quiescence and spikes characteristic of the recurrence inside the window. 

These solutions are particularly striking due to the wealth of complicated dynamics which can exist in models with more than two variables \cite{WIG2003}. The cause of the slow-fast recurrence-like behavior seems to overpower strange dynamics visible in $\delta=0.12$ and $\delta=0.1$, contrasting the multitude of peaks in the simulations of model \eqref{Eqn6} visible in Figure \ref{Fig9}. A comparison of these two models may be pivotal in understanding why the window seems to give rise to slow-fast dynamics.

\section{``Semi'' Recurrence} 

In this section, we present two systems that exhibit what we have called ``Semi''-recurrence. Systems with such semi-recurrence 
property do not have a sudden stop in continuous oscillations 
as does in the typical recurrence phenomenon, instead one will see 
weak slow-fast oscillating behaviour for certain values of 
the bifurcation parameter. What is interesting is that this ``weak'' behavior appears with a weak satisfaction of our hypotheses. 
The two models given in the following subsections show this weak
recurrence. As in the previous sections, we will present their bifurcation diagrams to show
a window which has periodic solutions bifurcating from a Hopf critical point followed by an additional Hopf bifurcation reducing the periodic solutions to damped oscillations. The damped oscillations are then stopped shortly by another bifurcation point. We are not interested in categorizing the bifurcation where the damped oscillations stop, simply remarking that they do.

\subsection{A controlled biochemical system model}

In this subsection, we study the following multi-component chemical control system with negative feedback which was studied in \cite{Walter1970} to explore the existence and ``spikey`` rhythms that exist in some biological systems. The dynamics of $n+1$ chemical species are given by
\begin{equation}\label{Eqn8} 
\begin{array}{rl} 
\dfrac{d S_1}{dt}+b_1S_1 \!\!\! &= \dfrac{K}{1+\alpha S_{n+1}^\rho}, \\[2.0ex]
\dfrac{d S_i}{dt}+b_iS_i \!\!\! &= b_{i-1}S_{i-1}, \quad i = 2, \dots, n+1.
\end{array}
\end{equation}
Each $b_i$ represents the first order rate constant at which species $S_i$ exits the systems for $i=1, \dots, n$. $K$ is the maximum uninhibited rate of the step inhibited by the final species $S_{n+1}$, while $\rho$ and $\alpha$ are the stoichiometry of the feedback inhibition and macroscopic inhibition, respectively.

All parameters take real positive values and all positive solutions stay positive. 
For our study we choose $n=5$, set 
$$
\alpha=1, \quad K=50, \quad \rho=9, \quad b_1=b_2=b_3=b_4=b_5=1,
$$ 
and treat $b_6$ as a bifurcation parameter. 
In order to find the equilibria of system \eqref{Eqn8}, we set $S_i'=0$, and then solve the set of equations: 
$$ 
b_1S^*_1 = \frac{K}{1+\alpha (S^*_6)^p}, \quad 
b_iS^*_i = b_{i-1}S^*_{i-1}, \quad i=2,\dots,6, 
$$ 
to obtain $S^*_{i+1}=\frac{b_i}{b_6}S_i^*$ for $i=1,\dots,5$. 
Substituting these solutions into the first equation in \eqref{Eqn8} yields the problem of solving the nonlinear equation for $S^*_6$:
$$ 
b_6S_6^* \big[1+\alpha (S_6^*)^p \big]=K.
$$
It can be shown that this equation has one positive and one negative solution for any values of $b_6>0$ and thus the system \eqref{Eqn8} has a unique positive equilibrium. When $b_6=0$, no equilibria exist. 

\begin{figure}[!h] 
\vspace*{0.5in} 
\begin{center}
\begin{overpic}[width=0.45\textwidth]{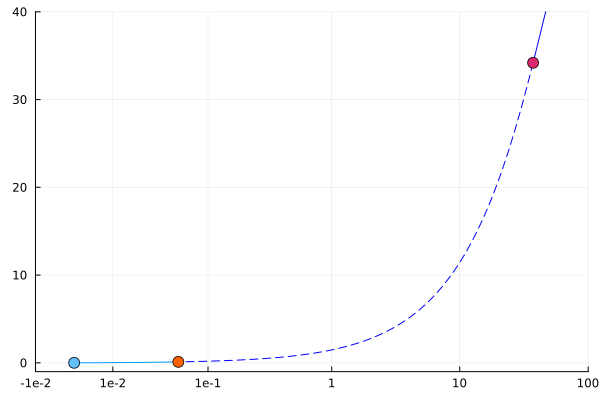} 
\put(5,70){$S_1$} 
\put(103,4){$b_6$} 
\end{overpic} 
\hspace*{0.20in} 
\begin{overpic}[width=0.45\textwidth]{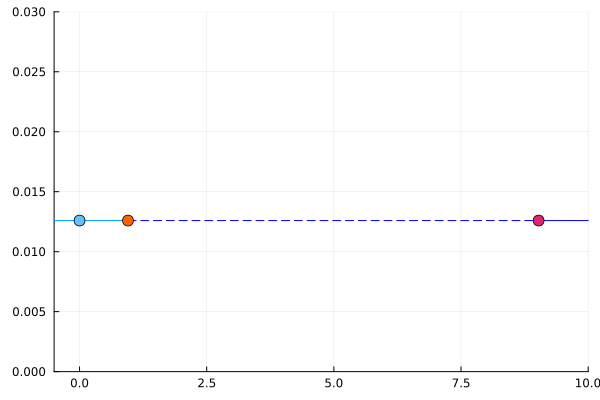}
\put(7,70){$y$} 
\put(103,4){$c_3$} 
\end{overpic} 

\hspace*{0.0in} (a) Bifurcation diagram of \eqref{Eqn8} \hspace*{0.10in} (b) Bifurcation diagram of \eqref{Eqn9}
\caption{Bifurcation diagrams of the models \eqref{Eqn8} and \eqref{Eqn9} showing a second Hopf bifurcation in close proximity to the stop point in the damped oscillations. The solid and dashed curves/lines indicate the stable and unstable equilibrium solutions, respectively. The bifurcation diagram in (a) is plotted on a symmetric logarithmic scale to better show the proximity of the stop in damped oscillations at $b_6=0$ (blue circle) to the Hopf bifurcation point (orange circle) as well as the second Hopf bifurcation (red circle). A similar explanation applies to (b). Note that if the small stable region between the Hopf bifurcations (red circles) and the stop in damped oscillations (blue circles) failed to exist, these models would satisfy our original hypotheses.}
         
\label{Fig11}
\end{center}

\end{figure}

\begin{figure}[h!] 
\vspace{0.20in} 
\begin{center} 
\begin{overpic}[width=0.28\textwidth]{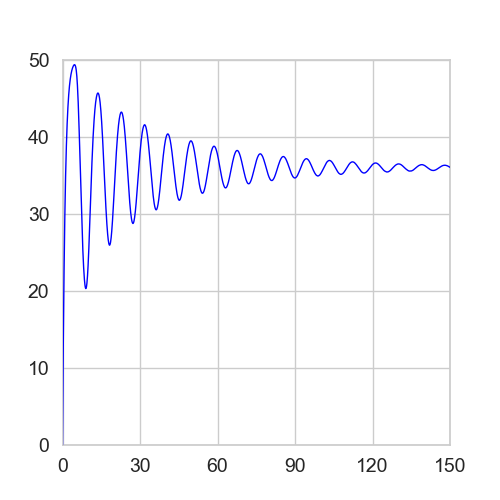} 
\put(8,95){$S_1$} 
\put(100,4){$t$} 
\end{overpic}
\hspace*{0.20in} 
\begin{overpic}[width=0.28\textwidth]{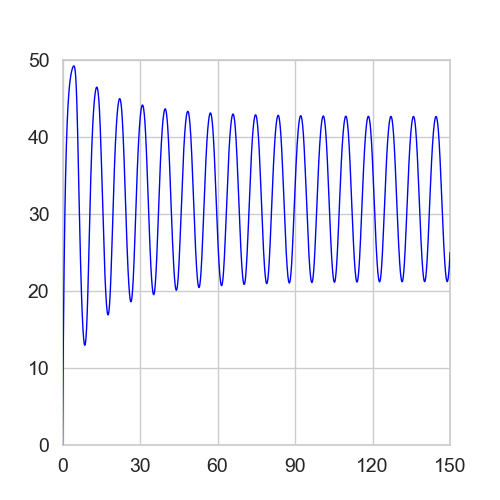}
\put(8,95){$S_1$} 
\put(100,4){$t$} 
\end{overpic} 
\hspace*{0.20in} 
\begin{overpic}[width=0.28\textwidth]{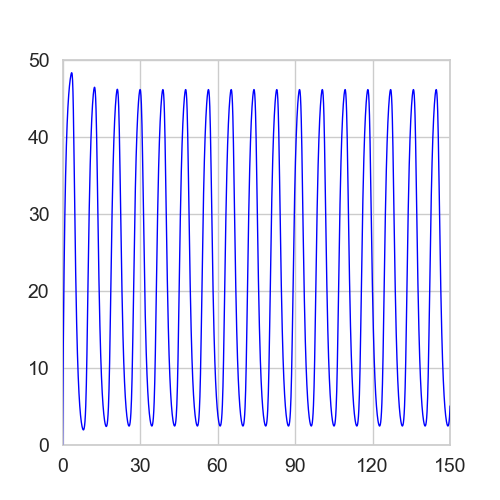}
\put(8,95){$S_1$} 
\put(100,4){$t$} 
\end{overpic} 

\hspace*{0.15in} (a) $b_6=40$ \hspace*{1.20in} (b) $b_6=35$
\hspace*{1.20in} (c) $b_6=20$

\vspace*{0.30in} 
\begin{overpic}[width=0.28\textwidth]{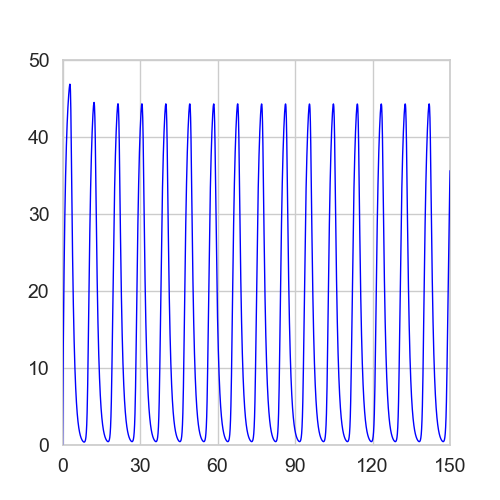}
\put(8,95){$S_1$} 
\put(100,4){$t$} 
\end{overpic} 
\hspace*{0.20in} 
\begin{overpic}[width=0.28\textwidth]{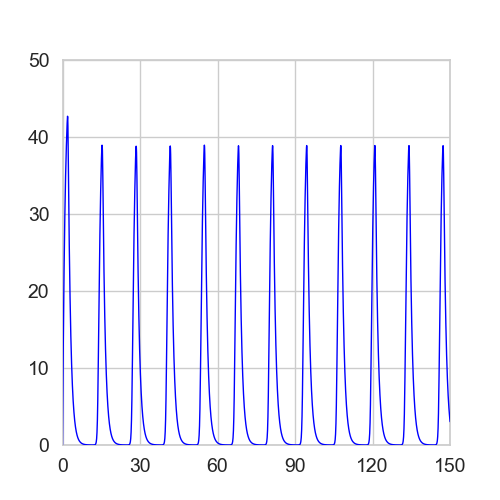}
\put(8,95){$S_1$} 
\put(100,4){$t$} 
\end{overpic} 
\hspace*{0.20in} 
\begin{overpic}[width=0.28\textwidth]{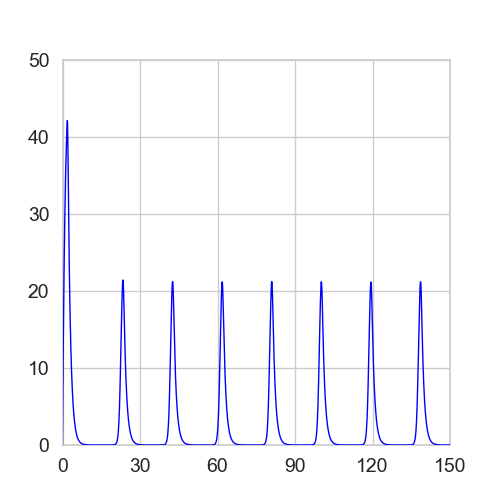}
\put(8,95){$S_1$} 
\put(100,4){$t$} 
\end{overpic} 

\hspace*{0.15in} (d) $b_6=10$ \hspace*{1.20in} (e) $b_6=1$
\hspace*{1.20in} (f) $b_6=0.3$

\vspace*{0.30in} 
\begin{overpic}[width=0.28\textwidth]{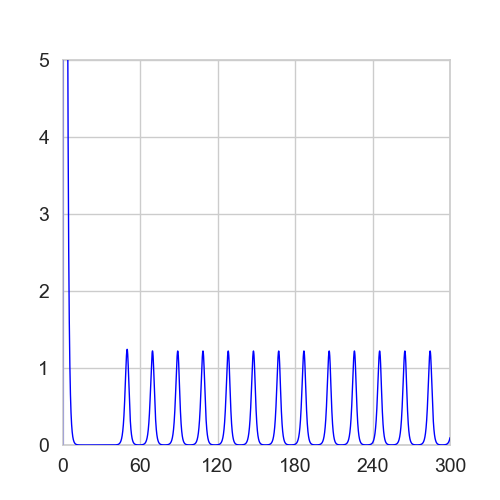}
\put(8,95){$S_1$} 
\put(100,4){$t$} 
\end{overpic} 
\hspace*{0.20in} 
\begin{overpic}[width=0.28\textwidth]{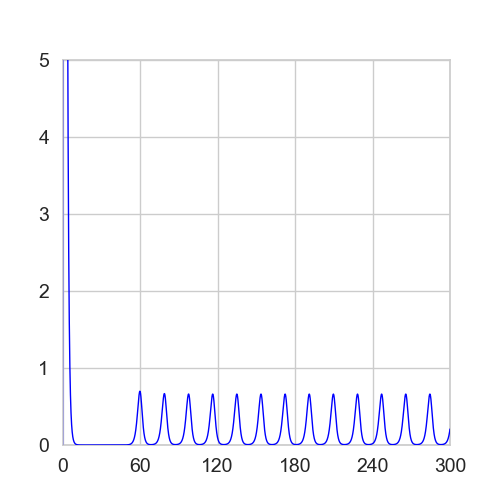}
\put(8,95){$S_1$} 
\put(100,4){$t$} 
\end{overpic} 
\hspace*{0.20in} 
\begin{overpic}[width=0.28\textwidth]{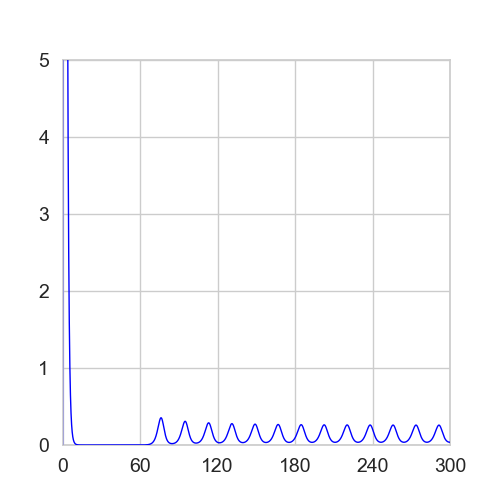}
\put(8,95){$S_1$} 
\put(100,4){$t$} 
\end{overpic} 

\hspace*{0.15in} (g) $b_6=0.1$ \hspace*{1.20in} (h) $b_6=0.08$
\hspace*{1.20in} (i) $b_6=0.06$

\vspace*{0.30in} 
\begin{overpic}[width=0.28\textwidth]{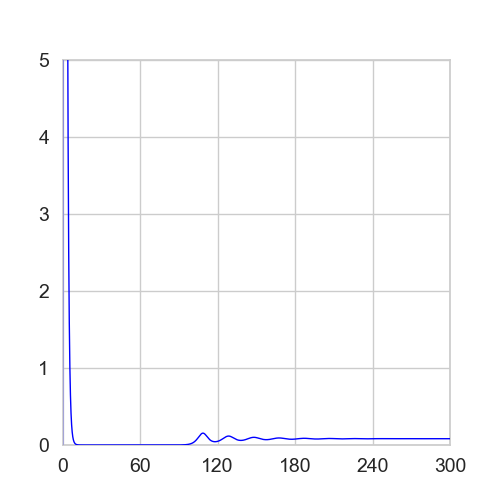}
\put(8,95){$S_1$} 
\put(100,4){$t$} 
\end{overpic} 
\hspace*{0.20in} 
\begin{overpic}[width=0.28\textwidth]{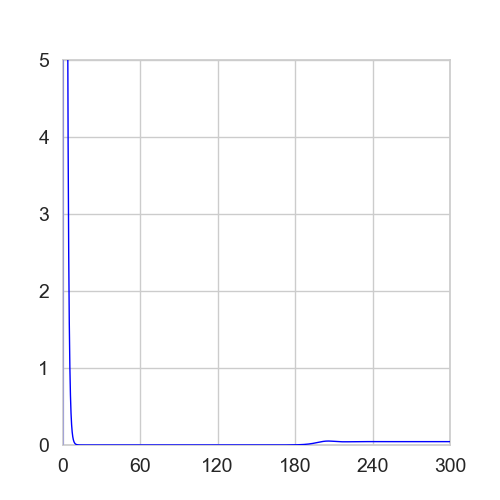}
\put(8,95){$S_1$} 
\put(100,4){$t$} 
\end{overpic} 
\hspace*{0.20in} 
\begin{overpic}[width=0.28\textwidth]{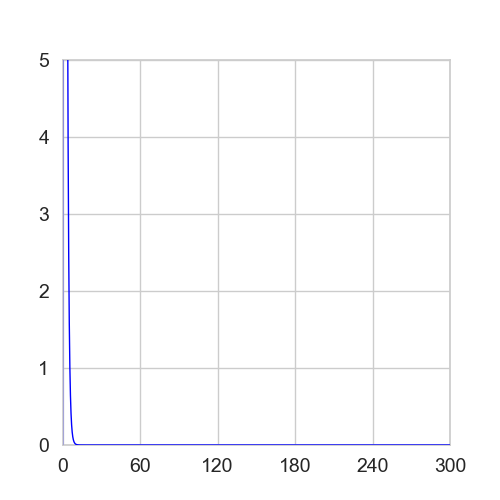}
\put(8,95){$S_1$} 
\put(100,4){$t$} 
\end{overpic} 

\hspace*{0.15in} (j) $b_6=0.04$ \hspace*{1.20in} (k) $b_6=0.02$
\hspace*{1.20in} (l) $b_6=0$

\caption{Simulated time histories of the model \eqref{Eqn8} for $\alpha=1$, $K=50$, $\rho=9$, $b_1=b_2=b_3=b_4=b_5=1$, with varying 
$b_6$ and the initial point: $S_i(0)=0$, $i=1, \dots, 6.$}
\label{Fig12}
\end{center} 

\end{figure}

As usual, we construct a bifurcation diagram numerically which is shown in figure \ref{Fig11}. The window condition is almost satisfied as we have an interval between the point $b_6=0$ and a Hopf bifurcation point at $b_6\approx 37.2$, where oscillations can exist, however a supercritical Hopf bifurcation existing at $b_6 \approx 0.05$ violates the window condition. Despite this, there is similar behaviour in the solutions as shown in Figure \ref{Fig12}, where the solutions display the spiked character of solutions in section 2 as $b_6$ approaches $0$. It seems to imply that the supercritical Hopf bifurcation at $b_6=0.05$ limits the length of the period of quiescence, causing solutions to die out slowly rather than suddenly. 

\subsection{ An SIR model} 

Finally, we consider an SIR model that demonstrates
secondary transmission routes to recurrent epidemics. The following is the reduced, dimensionless 
SIR model studied in \cite{Adams_Obara2011}: 
\begin{equation}\label{Eqn9} 
\begin{array}{ll}
\dfrac{dx}{dt} =(x+y)(1-x-y)-c_1xy-xz, \\[1.5ex] 
\dfrac{dy}{dt} = c_1xy+xz-c_2y, \\[1.5ex] 
\dfrac{dz}{dt} = c_3(c_4y-z),
\end{array}
\end{equation}
where $c_1$, $c_2$, $c_3$ and $c_4$ are positive parameters. It has been shown in \cite{Adams_Obara2011} that the model has four equilibria, $(0,0,0)$, $(1,0,0)$,
$E_+=(x_0, y_{0+}, c_4y_{0+})$, and $E_- =(x_0, y_{0-}, c_4y_{0-})$, where 
$$
x_0=\frac{c_2}{c_1+c_4} \quad 
\textrm{and} \quad 
y_{0\pm} = \frac{1}{2}(1-2x_0-c_2)\pm \frac{1}{2} \sqrt{(2x_0 + c_2-1)^2-4x_0(x_0)-1}.
$$
We focus on the parameter values for which
$E_-$ is irrelevant. In \cite{Adams_Obara2011} the author showed that the limit cycles bifurcate from the equilibrium $E_+$ due to a Hopf bifurcation. 
For simplicity, we fix
\[c_1=1, \ \ c_2=5, \ \ c_4=90, \]
and treat $c_3$ as a bifurcation parameter.
As we did before, a bifurcation diagram is numerically constructed and presented in Figure \ref{Fig11}(b). Similarly to \ref{Fig11}(a), two Hopf bifurcations enclose a region where oscillations can exist.

Like discussed above, there seems to be a limit to how ``spikey'' the solution can become as the parameter approaches the subcritical Hopf bifurcation point. Looking closer at the bifurcation diagram one can observe three bifurcation points: 
A Hopf bifurcation at $c_3=9.0259$ generating a limit cycle, a second Hopf bifurcation at $c_3=0.9523$ destroying the limit cycle, and then a bifurcation point where the decaying oscillations cease at $c_3=0$.
\begin{figure}[h!] 
\vspace{0.20in} 
\begin{center} 
\begin{overpic}[width=0.28\textwidth]{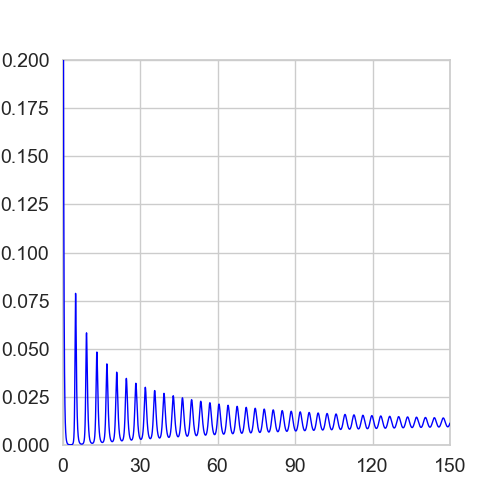} 
\put(10,95){$y$} 
\put(100,5){$t$} 
\end{overpic}
\hspace*{0.20in} 
\begin{overpic}[width=0.28\textwidth]{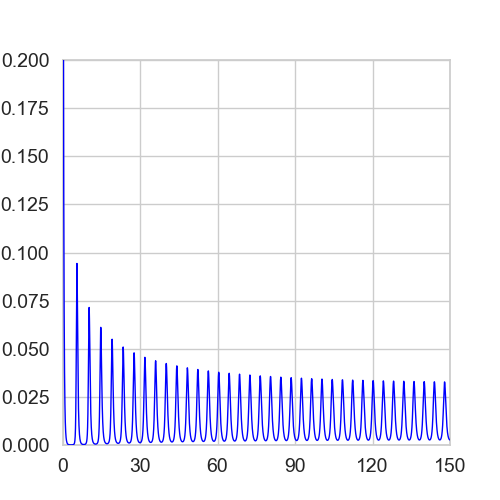}
\put(10,95){$y$} 
\put(100,5){$t$} 
\end{overpic} 
\hspace*{0.20in} 
\begin{overpic}[width=0.28\textwidth]{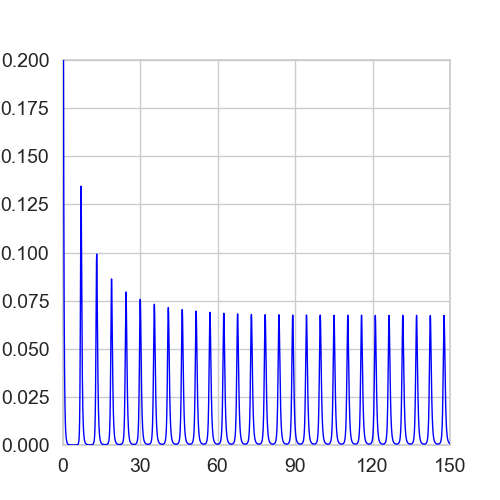}
\put(10,95){$y$} 
\put(100,5){$t$} 
\end{overpic} 

\hspace*{0.15in} (a) $c_3=12$ \hspace*{1.20in} (b) $c_3=8$
\hspace*{1.20in} (c) $c_3=4$

\vspace*{0.30in} 
\begin{overpic}[width=0.28\textwidth]{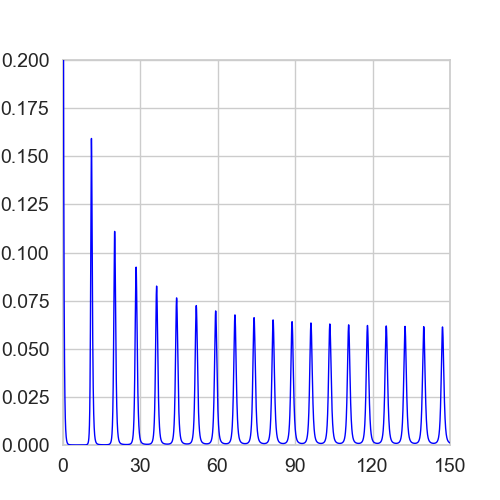} 
\put(10,95){$y$} 
\put(100,5){$t$} 
\end{overpic} 
\hspace*{0.20in} 
\begin{overpic}[width=0.28\textwidth]{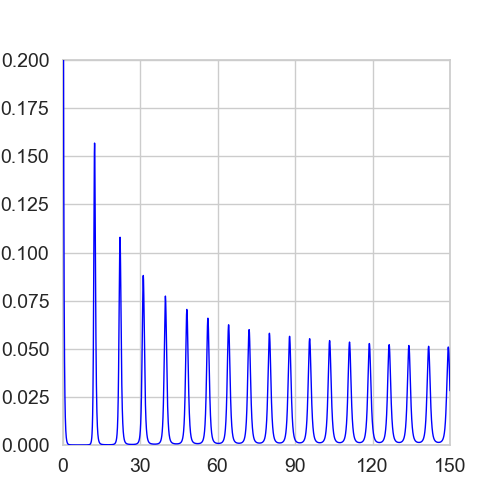} 
\put(10,95){$y$} 
\put(100,5){$t$} 
\end{overpic} 
\hspace*{0.20in} 
\begin{overpic}[width=0.28\textwidth]{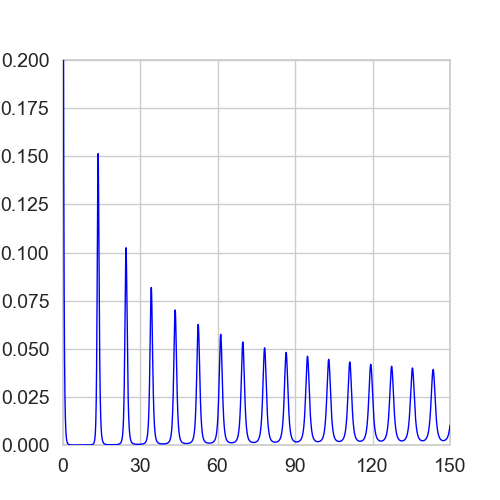} 
\put(10,95){$y$} 
\put(100,5){$t$} 
\end{overpic} 

\hspace*{0.15in} (d) $c_3=1.5$ \hspace*{1.20in} (e) $c_3=1.2$
\hspace*{1.20in} (f) $c_3=0.95$

\vspace*{0.30in} 
\begin{overpic}[width=0.28\textwidth]{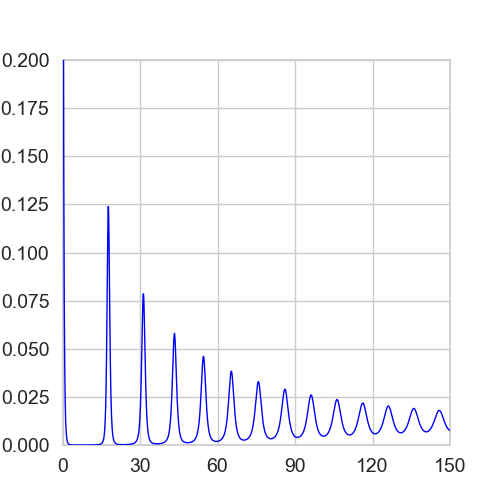} 
\put(10,95){$y$} 
\put(100,5){$t$} 
\end{overpic} 
\hspace*{0.20in} 
\begin{overpic}[width=0.28\textwidth]{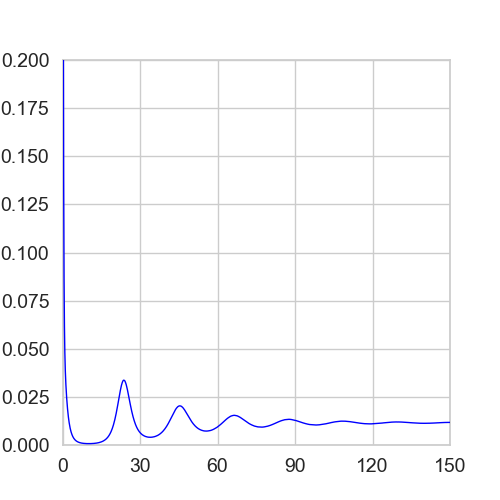} 
\put(10,95){$y$} 
\put(100,5){$t$} 
\end{overpic} 
\hspace*{0.20in} 
\begin{overpic}[width=0.28\textwidth]{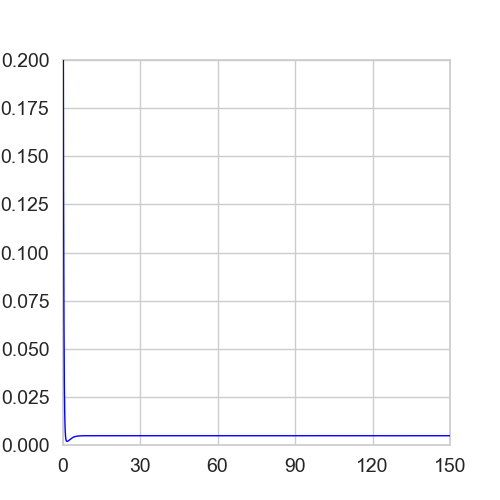} 
\put(10,95){$y$} 
\put(100,5){$t$} 
\end{overpic} 

\hspace*{0.15in} (g) $c_3=0.5$ \hspace*{1.20in} (h) $c_3=0.1$
\hspace*{1.20in} (i) $c_3=0$

\caption{Simulated time history of the model \ref{Eqn9} 
for $c_1=1$, $c_2=5$, $c_3=90$ with varying $c_3$ 
and the initial point $(x(0),y(0),z(0))=(0.09777,0.02081,1.04290)$.}
\label{fig:first}
\end{center} 

\end{figure}

\section{Conclusion and Discussion}

Proving the existence of recurrence in real world models can be a complicated task, requiring the use of GSPT or extending it as in \cite{Hsu2019, Li2016}. The simple criterion developed in \cite{ZWY2014,YW2019,YZ2019} based on bifurcation theory can be easily applied to identify slow-fast motions in practical systems. In this paper, we present five models that fail to satisfy this criterion but still
show recurrence. Therefore, we generalize those hypotheses so that they can be used to analyze more systems like those five models for the existence of recurrence. In addition, we showed semi-recurrence in two examples with windows to almost match our new hypotheses, 
suggesting that the recurrence phenomenon may be the limit
of a weaker condition. Although we have provided the evidence that
the recurrence phenomenon is dependent on these conditions, we have yet to establish a satisfactory mathematical theory for our hypotheses.

The examples we have presented in this paper, as well as those discussed in \cite{YZ2019}, allow us to look for heuristic arguments that could motivate developing a theory in the future. 
Both the conditions given in \cite{YZ2019} and the hypotheses proposed in this paper focus on the bifurcation of a single equilibrium solution. 
The condition about the window to generate a continuous oscillation
implies that the Hopf bifurcation results in a birth to producing limit cycles. In addition, we imposed the condition that the bifurcating limit cycles must disappear ``suddenly'', contrasting a gradual disappearance characteristic of a subcritical Hopf bifurcation. In \cite{YZ2019}, the authors provided evidence to support that such a phenomenon should be due to the equilibrium undergoing a transcritical or saddle-node bifurcation. 
However, we suggest that the transcritical or saddle-node bifuration condition can be better explained by interactions of the limit cycle with nearby stable and unstable manifolds.

\begin{figure}[!h] 
\vspace*{0.2in} 
\begin{center}
\begin{overpic}[width=0.45\textwidth]{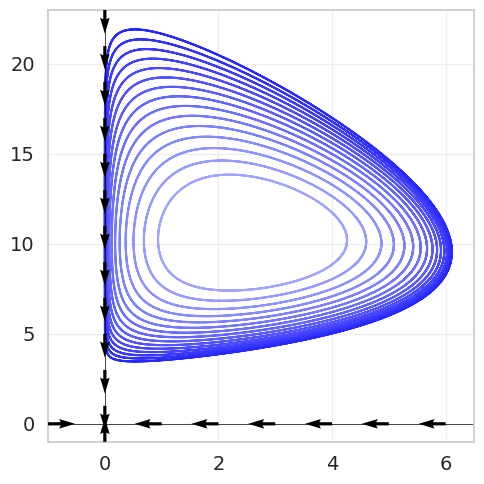} 
\put(7,100){$y$} 
\put(100,5){$x$} 
\end{overpic} 
\hspace*{0.20in} 
\begin{overpic}[width=0.45\textwidth]{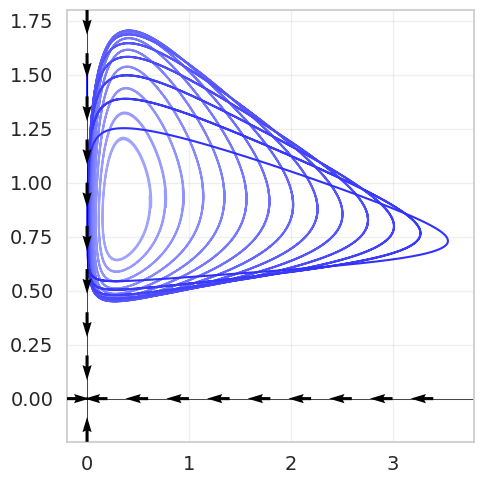} 
\put(9,100){$y$} 
\put(100,5){$x$} 
\end{overpic} 

{
\hspace*{0.20in}(a) Limit cycles of model \eqref{Eqn3} (see Figure \ref{Fig1}) \hspace*{0.45in} (b) Limit cycles of model \eqref{Eqn5} (see Figure \ref{Fig5}).
        }
\caption{ Limit cycles bifurcating in models \eqref{Eqn3} and \eqref{Eqn5} as the bifurcation parameter $\varepsilon$ is varied within the window. The stable and unstable manifolds of the equilibrium at $(0,0)$ are displayed as arrows along the axes. Periods of quiescence coincide with the points in which the limit cycles come in close contact to the stable/unstable manifolds. As the bifurcation parameter approaches the stop in oscillations, the limit cycle grows to be in closer and closer contact with the stable/unstable manifolds. The hypotheses of \cite{YZ2019} tended to work because the transcritical and saddle-node bifurcations force limit cycles to be in close contact with stable/unstable manifolds of another equilibrium.} 
         
\label{Fig16}
\end{center} 

\end{figure}

In \cite{ZWY2013,ZWY2014,YW2019,YZ2019}, the authors used bifurcation theory, based on eigenvalues, to explain why a transcritical or a saddle-node bifurcation is necessary for the appearance of recurrence, which leaves unexplained recurrences where transcritical and saddle-node bifurcations are absent. To rectify this, we suggest that the recurrence is due to stable and unstable manifolds being in close proximity to the limit cycle produced by the Hopf bifurcation in our theory. To illustrate this, the limit cycles of systems \eqref{Eqn3} and \eqref{Eqn5} are graphed as the bifurcation parameter is varied within the window in Figure \ref{Fig16}. As one can see, the limit cycles seem to have some tendency to grow but are limited by the nearby stable/unstable manifold of the equilibria at zero due to uniqueness of solutions. Similar behavior can be seen in the other systems discussed in our article. This suggests that a theory capturing this phenomenon should focus on the interaction between the stable manifold and the limit cycles emerging from the Hopf bifurcation. More work is needed to develop a fundamental mathematical theory in the near future. 

\section*{Acknowledgement}  
This work was supported by the Natural Sciences and Engineering Research Council of Canada (No.~R2686A02).

%
%



\begin{thebibliography}{00}

\bibitem{Adams_Obara2011}
M.~R. Adams, S. Obara,
Recurrent epidemics resulting from secondary transmission routes in SIR models,
Electron. J. Diff. Equ. 2011(142) (2011) 1-12.

\bibitem{ALL1996}
Alligood, K.T., Sauer, T.D., Yorke, J.A., 1996. Chaos: An Introduction to Dynamical Systems, Textbooks in Mathematical Sciences. Springer, New York, NY. https://doi.org/10.1007/b97589

\bibitem{ARN1989}
Arnold, V.I. (Vladimir I., 1989. Mathematical methods of classical mechanics. New York: Springer-Verlag.

\bibitem{CAR1981}
Carr, J., 1981. Applications of Centre Manifold Theory, Applied Mathematical Sciences. Springer US, New York, NY. https://doi.org/10.1007/978-1-4612-5929-9

\bibitem{DIM2003}
Di Mascio, M., Markowitz, M., Louie, M., Hogan, C., Hurley, A., Chung, C., Ho, D.D., Perelson, A.S., 2003. Viral Blip Dynamics during Highly Active Antiretroviral Therapy. J Virol 77, 12165–12172. https://doi.org/10.1128/JVI.77.22.12165-12172.2003

\bibitem{EDE2020}
Edelstein-Keshet, L., 2020. Differential Calculus for the Life Sciences.

\bibitem{ERD1989}
Érdi, P., Tóth, J., 1989. Mathematical Models of Chemical Reactions: Theory and Applications of Deterministic and Stochastic Models. Manchester University Press.

\bibitem{Fenichel1979} 
N. Fenichel, 
Geometric singular perturbation theory for ordinary differential equations,
J. Differ. Equ. 31(1) (1979) 53-98. 
doi: 10.1016/0022-0396(79)90152-9.  

\bibitem{GOO2011}
Goodwine, B., 2011. Engineering Differential Equations: Theory and Applications. Springer, New York, NY. https://doi.org/10.1007/978-1-4419-7919-3

\bibitem{GUC1983}
Guckenheimer, J., Holmes, P., 1983. Nonlinear Oscillations, Dynamical Systems, and Bifurcations of Vector Fields, Applied Mathematical Sciences. Springer, New York, NY. https://doi.org/10.1007/978-1-4612-1140-2

\bibitem{HEK2010}
Hek, G., 2010. Geometric singular perturbation theory in biological practice. J. Math. Biol. 60, 347–386. https://doi.org/10.1007/s00285-009-0266-7

\bibitem{Hsu2019}
T.-H. Hsu,
Number and stability of relaxation oscillations for predator-prey 
systems with small death rates,
SIAM J. Appl. Math. 18(1) (2019) 33-67.
doi: 10.1137/18M1166705.

\bibitem{Hsu2015}
S.-B. Hsu, S. Ruan, T.-H. Yang,
Analysis of three species Lotka–Volterra food web models with omnivory,
J. Math. Anal. Appl. 426(2) (2015) 659-687. 

\bibitem{JAC2005}
Jacobson, M.Z., 2005. Fundamentals of Atmospheric Modeling, 2nd ed. Cambridge University Press, Cambridge. https://doi.org/10.1017/CBO9781139165389

\bibitem{JON1995}
Jones, C.K.R.T., 1995. Geometric singular perturbation theory, in: Arnold, L., Jones, C.K.R.T., Mischaikow, K., Raugel, G., Johnson, R. (Eds.), Dynamical Systems: Lectures Given at the 2nd Session of the Centro Internazionale Matematico Estivo (C.I.M.E.) Held in Montecatini Terme, Italy, June 13–22, 1994. Springer, Berlin, Heidelberg, pp. 44–118. https://doi.org/10.1007/BFb0095239

\bibitem{JUL2017}
Bezanson, J., Edelman, A., Karpinski, S., Shah, V.B., 2017. Julia: A fresh approach to numerical computing. SIAM Review 59, 65–98. https://doi.org/10.1137/141000671

\bibitem{Kaklamanos2023}
P. Kaklamanos, N. Popovi\'{e},
Complex oscillatory dynamics in a three-timescale El Ni$\tilde{\rm n}$o 
Southern Oscillation model, 
Physica D 449 (2023) 133740. 
doi: 10.1016/j.physd.2023.133740. 

\bibitem{KAR2012}
Karr, J.R., Sanghvi, J.C., Macklin, D.N., Gutschow, M.V., Jacobs, J.M., Bolival, B., Assad-Garcia, N., Glass, J.I., Covert, M.W., 2012. A Whole-Cell Computational Model Predicts Phenotype from Genotype. Cell 150, 389–401. https://doi.org/10.1016/j.cell.2012.05.044

\bibitem{KLI2016}
Klipp, E., Liebermeister, W., Wierling, C., Kowald, A., 2016. Systems Biology: A Textbook. John Wiley \& Sons.

\bibitem{KUZ2023}
Kuznetsov, Y.A. (2023). Numerical Analysis of Bifurcations. In: Elements of Applied Bifurcation Theory. Applied Mathematical Sciences, vol 112. Springer, Cham. https://doi.org/10.1007/978-3-031-22007-4\_10

\bibitem{LI2013}
Li, C., Zhu, H., 2013. Canard cycles for predator–prey systems with Holling types of functional response. Journal of Differential Equations 254, 879–910. https://doi.org/10.1016/j.jde.2012.10.003

\bibitem{LI2016}
Li, C., 2016. Slow-Fast Dynamics and Its Application to a Biological Model, in: Toni, B. (Ed.), Mathematical Sciences with Multidisciplinary Applications. Springer International Publishing, Cham, pp. 301–325. https://doi.org/10.1007/978-3-319-31323-8\_14

\bibitem{Li2016}
M.~Y. Li, W. Liu, C. Shan, Y. Yi,
Turning points and relaxation oscillation cycles in simple epidemic models,
SIAM J. Appl. Math. 76(2) (2016) 663-687.
doi: 10.1137/15M1038785.

\bibitem{LOR1963}
Lorenz, E.N., 1963. Deterministic Nonperiodic Flow. Journal of the Atmospheric Sciences 20, 130–141. https://doi.org/10.1175/1520-0469(1963)020<0130:DNF>2.0.CO;2

\bibitem{MIS1980}
Mishchenko, E.F., Rozov, N.Kh., 1980. Differential Equations with Small Parameters and Relaxation Oscillations. Springer US, Boston, MA. https://doi.org/10.1007/978-1-4615-9047-7

\bibitem{MUR2002}
Murray, J.D. (Ed.), 2002. Mathematical Biology: I. An Introduction, Interdisciplinary Applied Mathematics. Springer, New York, NY. https://doi.org/10.1007/b98868

\bibitem{Sahoo-Samanta2023}  
D. Sahoo, G. Samanta,
Oscillatory and transient dynamics of a slow–fast predator–prey system 
with fear and its carry-over effect, 
Nonlinear Analysis: Real World Appl. 73 (2023) 103888.
doi: 10.1016/j.jmaa.2015.01.035. 

\bibitem{SCI2020}
Virtanen, P., Gommers, R., et al., 2020. SciPy 1.0: fundamental algorithms for scientific computing in Python. Nat Methods 17, 261–272. https://doi.org/10.1038/s41592-019-0686-2

\bibitem{SOR2016}
Sörstedt, E., Nilsson, S., Blaxhult, A., Gisslén, M., Flamholc, L., Sönnerborg, A., Yilmaz, A., 2016. Viral blips during suppressive antiretroviral treatment are associated with high baseline HIV-1 RNA levels. BMC Infect Dis 16, 305. https://doi.org/10.1186/s12879-016-1628-6

\bibitem{SUN2022}
Sundararajan, Dr.D., 2022. Control Systems: An Introduction. Springer International Publishing, Cham. https://doi.org/10.1007/978-3-030-98445-8


\bibitem{BIF2020}
Veltz, R., 
BifurcationKit.jl.
2020. 

\bibitem{Walter1970}  
C.~F. Walter,
The occurrence and the significance of limit cycle behavior in controlled biochemical systems, 
J. Theor. Biol. 27(2) (1970) 259-272. 
Doi: 10.1016/0022-5193(70)90141-4. 

\bibitem{WIG2003}
Wiggins, S., 2003. Introduction to Applied Nonlinear Dynamical Systems and Chaos, Texts in Applied Mathematics. Springer-Verlag, New York. https://doi.org/10.1007/b97481

\bibitem{YW2019}
P. Yu, X. Wang, 
Analysis on recurrence behavior in oscillating networks of biologically 
relevant organic reactions,
Math. Biosci. Eng. 16(5) (2019) 5263-5286. 
doi: 10.3934/mbe.2019263. 

\bibitem{YZ2019}
P. Yu, W. Zhang, 
Complex dynamics in a unified SIR and HIV disease model: 
A bifurcation theory approach,
J. Nonlinear Sci. 29(5) (2019) 2447-2500.
doi:10.1007/s00332-019-09550-7. 


\bibitem{ZWY2013} 
W. Zhang, L. Wahl, P. Yu, 
Conditions for transient viremia in deterministic in-host models: 
viral blips need no exogenous trigger,
SIAM J. Appl. Math. 73(2) (2013) 853-881.
doi:10.2307/23479953. 

\bibitem{ZWY2014} 
W. Zhang, L. Wahl, P. Yu, 
Viral blips may not need a trigger: How transient viremia can 
arise in deterministic in-host models,
SIAM Review 56(1) (2014) 127-155. 
doi:10.1137/130937421.
\end{thebibliography}



\end{document}